\documentclass[winedt,yap]{iitparc}
\usepackage{amssymb}
\usepackage{graphics}
\usepackage{graphicx}

\begin{document}

\frontmatter          

\IssuePrice{25}%
\TransYearOfIssue{2001}%
\TransCopyrightYear{2001}%
\OrigYearOfIssue{2001}%
\OrigCopyrightYear{2001}%
\TransVolumeNo{62}%
\TransIssueNo{3}%
\OrigIssueNo{3}%
\OrigPages{108--133}%

\mainmatter              


\Rubrika{EVOLVING SYSTEMS} \CRubrika{EVOLVING SYSTEMS}

\setcounter{page}{443}

\newcommand{\PLE}[1]{{\bf Proof of Lemma~#1.}}

\def\epr{\hfill$\sqcap\!\!\!\!\sqcup$\medskip\par}

\newcommand\propro[2]{\medskip\par\PPR{#1}{\rm #2}\par}
\newcommand\prothe[2]{\medskip\par\PTH{#1}{\rm #2}\par}
\newcommand\prolem[2]{\medskip\par\PLE{#1}{\rm ~#2}\par}
\newcommand\theoremprimed[2]{\par\hskip5mm{\bfseries Theorem #1.\ \ }{\it #2}\par\bigskip}
\newcommand\coroll[2]{\par{\bfseries Corollary #1.\ \ }{\itshape #2}\par\bigskip}
\newcommand\axiom[2] {\bigskip\par{\bfseries #1.\ } {\rm #2}\par\bigskip}
\newcommand\axiomt[2]{\par{\bfseries #1.\ }         {\rm #2}\par\bigskip}

\def\statemento#1#2{\par{\sl #1 }{\it #2}\par}    
\def\secondstatemento#1#2{{\sl #1 }{\it #2}\par}  
\def\mo{\Big\|}                                   
\def\moo{\vspace{-0.9ex}\Big\|}                   
\def\D{{\rm\Delta}}                               
\def\l{\ell}                                      
\def\aa{\alpha}                                   
\def\G{\Gamma}                                    
\def\la{\lambda}                                  
\def\si{\sigma}                                   
\def\e{\varepsilon}                               
\def\h{H}                                         
\def\FF{\mathop{\cal F}\nolimits}                 
\def\PP{\mathop{\cal P}\nolimits}                 
\def\GG{\mathop{\cal G}\nolimits}                 
\def\TT{\mathop{\cal T}\nolimits}                 
\def\NN{\mathop{\cal N}\nolimits}                 
\def\RR{\mathop{\cal R}\nolimits}                 
\def\ul   {\mathop{\rm u}  \limits}               
\def\suml {\mathop{\sum}   \limits}               
\def\liml {\mathop{\lim}   \limits}               
\def\capo {\mathop{\bigcap}\limits}               
\def\cupo {\mathop{\bigcup}\limits}               
\def\maxl {\mathop{\max}   \limits}               
\def\ulo{\ul^{\scriptscriptstyle \circ}}          
\def\cdc{,\ldots,}                                
\def\1n{1,\ldots,n}                               
\def\_#1{\mathop{\hspace{-2pt}^{}_{#1}}}          
\def\epr{\hfill$\square$\smallskip\par}           
\def\parr{\newline\indent}                        
\def\R{{\mathbb R}}                               
\def\C{{\mathbb C}}                               
\def\oj {\mathop{\bar{J}}\nolimits}               
\def\vj {\mathop{\bar{J}}\nolimits}               
\def\q  {\bar{J}}                                 
\def\ktil {\widetilde K}                          
\def\Ptil {\widetilde P}                          
\def\ptil {\widetilde p}                          
\def\ltil {L^{\scriptscriptstyle (+)}}            
\def\qtil {\widetilde Q}                          
\def\intercal{\mathop{\scriptscriptstyle T}\nolimits}
\def\Apr{A^{\scriptscriptstyle\#}}                
\def\Lpr{L^{\scriptscriptstyle\#}}                
\def\lpr{\l^{\scriptscriptstyle\#}}               
\def\tr{{\rm tr}}                                 
\def\x{{}}                                        
\def\z{{}}                                        
\def\NL{{\cal N}({\bf L})}                        
\def\NLT{{\cal N}({\bf L}^{\intercal})}           
\def\NJ{{\cal N}({\bf\vj})}                       
\def\NJT{{\cal N}({\bf\vj}^{\intercal})}          
\def\RL{{\cal R}({\bf L})}                        
\def\RLT{{\cal R}({\bf L}^{\intercal})}           
\def\RJ{{\cal R}({\bf\vj})}                       
\def\RJT{{\cal R}({\bf\vj}^{\intercal})}
\raggedbottom

\author{R. P. Agaev \and P. Yu. Chebotarev}
\title{Spanning Forests of a Digraph\\
and Their Applications%
\thanks{This work was supported by the Russian
Foundation for Basic Research and INTAS. }}

\institute{Trapeznikov Institute of Control Sciences, Russian Academy of Sciences, Moscow, Russia
\\e-mail$:$ pchv@rambler.ru$,$ chv@lpi.ru
}

\received{Received September 18, 2000}
\titlerunning{Spanning Forests of a Digraph}
\authorrunning{Agaev, Chebotarev}
\OrigCopyrightedAuthors{Agaev, Chebotarev}
\OrigPages{108--133}%

\maketitle

\begin{abstract}
{We study spanning diverging forests of a digraph and related matrices. It is shown that the normalized
matrix of out forests of a digraph coincides with the transition matrix in a specific observation model for
Markov chains related to the digraph. Expressions{\x} are given for the Moore-Penrose generalized inverse and
the group inverse of the Kirchhoff matrix. These expressions involve the matrix of maximum out forests{\x} of
the digraph. Every matrix of out forests with a fixed number of arcs and the normalized matrix of out forests
are represented as polynomials of the Kirchhoff matrix; with the help of these identities,{\x} new proofs are
given for the matrix-forest theorem and some other statements. A connection is specified between the forest
dimension of a digraph and the degree of an annihilating polynomial for the Kirchhoff matrix. Some
accessibility measures for digraph vertices are considered. These are based on the enumeration of spanning
forests.

}
\end{abstract}

\section{Introduction}

Directed graphs provide a simple and universal tool to model connection structures. It is not accidental that
the first systematic monograph in the theory of digraphs~\cite{HarNoCa} was titled ``Structural Models: An
Introduction to the Theory of Directed Graphs.'' Digraphs frequently serve to model processes that can
proceed in the direction of arcs. Physical transference, service, control, transmission of influences, ideas,
innovations, and diseases are examples of such processes. If a process can start from a number of vertices
and ends with the inclusion of all vertices, then the process can be modelled by the family of out forests
(i.e., spanning diverging forests) of the digraph. The enumeration of all out forests allows one to determine
the typical roles of the vertices in the process: one vertex is a typical starting point, another vertex is a
typical intermediate point, some vertex is a typical terminating point of the process, etc. If an initial
(weighted) digraph imposes some measure on the said processes, then the ``role profile'' of each vertex can
be expressed numerically. Moreover, an exact answer can be given to the following important question: how
likely is it that the process initiated at vertex $j$ arrives at vertex~$i$. It is not surprising that out
forests of a digraph turn out to be closely related with Markov chains realizable on the digraph.

The study of out forests has been started in~\cite{Fiedler}. Generally, they were given less attention in the
literature, than that given to spanning diverging trees (out arborescences), which exist only for a narrow
class of digraphs.{\x} We mention in this connection
\cite{Kelm,LiuChow,Cha,MyrvoldA,BapatConstantine,Takacs,Erdos,Merr97,Merr98}, where still undirected forests
were considered in most cases. The{\x} maximum out forests (i.e., out forests with the greatest possible
number of arcs) of a digraph were studied in~\cite{Che.rat,Che.ra1}. It was established that the normalized
matrix of such forests coincides with the matrix of limiting probabilities of every Markov chain {\em
related\/} to the given digraph. Some results on spanning forests of directed and undirected multigraphs were
given in~\cite{CheSha97,CheSha981}.

In this paper, we study the normalized matrix of out forests (which{\x} has been also termed the matrix of
relative forest accessibilities and the matrix of forest proximities) and the matrices of forests with fixed
numbers of arcs.

\section{Notation and some earlier results}
\label{sec2}

In the terminology, we mainly follow~\cite{HarNoCa,Harary}. Suppose that $\G$ is a weighted digraph without
loops, $V(\G)=\{\1n\}$ $(n>1)$ is its set of vertices, and $E(\G)$ its set of arcs. The weights of all arcs
are supposed to be strictly positive. {\it A subgraph\/} of a digraph $\G$ is a digraph whose vertices and
arcs belong to the sets of vertices and arcs of $\G$;{\x} the weights of subgraph's arcs are the same as
in~$\G$. {\it A restriction\/} of $\G$ to $V'\subset V(\G)$ is a digraph whose arc set contains all the arcs
in $E(\G)$ that have both incident vertices in~$V'$. {\it A spanning subgraph\/} of $\G$ is a subgraph with
vertex set $V(\G)$. The {\it indegree\/} id($w$) of vertex $w$ is the number of arcs that come to~$w$, {\it
outdegree\/} od($w$) of vertex $w$ is the number of arcs that come from~$w$. A vertex $w$ will be called {\it
undominated\/} if id($w$)=0 and {\it dominated\/} if id$(w)\ge 1$. A vertex $w$ is {\it isolated\/} if~$\G$
contains no arcs incident to~$w$.

A {\it route\/} in a digraph is an alternating sequence of vertices and arcs $w_0,e_1,$ $w_1\cdc e_k, w_k$
with every arc $e_i$ being $(w_{i-1},w_i)$. If every arc $e_i$ is either $(w_{i-1},w_i)$ or $(w_i,w_{i-1}),$
then the sequence is called a {\it semiroute}. A {\it path\/} in a digraph is a route all whose vertices are
different. A~{\it circuit\/} is a route with $w\_0=w_k$, the other vertices being distinct and different from
$w_0$. A~vertex $w$ {\it is reachable\/} from a vertex $z$ in $\G$ if $w=z$ or $\G$ contains a path from $z$
to~$w$. A~{\it semicircuit\/} is an alternating sequence of distinct vertices and arcs, $w_0,e_1,w_1\cdc
e_k,w_0,$ where every arc $e_i$ is either $(w_{i-1},w_i)$ or $(w_i,w_{i-1}$) and all vertices $w_0\cdc
w_{k-1}$ are different. The restriction of $\G$ to any maximal subset of vertices connected by semiroutes is
called a {\it weak component\/} of~$\G$.{\x} Let $E=(\e_{ij})$ be the matrix of arc weights. Its entry
$\e_{ij}$ is zero if and only if there is no arc from vertex $i$ to vertex~$j$ in~$\G$. If $\G'$ is a
subgraph of $\G$, then the weight of $\G'$, $\e(\G')$, is the product of the weights of all its arcs; if
$\G'$ {\x}does not contain arcs, then $\e(\G')=1$.  The weight of a nonempty set of digraphs $\GG$ is defined
as follows:
\[
\e(\GG)=\suml_{H\in\GG}\e(H);
\]
the weight of the empty set is~0.

The {\it Kirchhoff matrix\/} \cite{Tutte} of a weighted digraph $\G$ is the $n\times n$-matrix
$L=L(\G)=(\l\_{ij})$ with elements $\l\_{ij}=-\e_{ji}$ when $j\ne i$ and $\l\_{ii}=-\suml_{k\ne i}\l\_{ik}$,
$i,j=\1n$.

A {\it diverging tree\/} is a digraph without semicircuits that has a vertex (called the {\it root}\/){\x}
from which every vertex is reachable. The indegree of every non-root vertex of a diverging tree is~1. If $w$
is the root, then id$(w)=0$. A {\it converging tree\/} is a digraph without semicircuits that has a vertex
(called the {\it sink}\/){\x} reachable from every {\x}vertex.

A~{\it diverging forest\/} ({\it converging forest}\/){\x} is a digraph without circuits such that
id$(w)\le1$ (respectively, od$(w)\le1$) for every {\x}vertex~$w$. An {\em out forest\/} ({\em in
forest}\/){\x} of a digraph $\G$ is any its spanning diverging (respectively, converging) forest.

The weak components of diverging forests (converging forests) are diverging trees (respectively,{\x}
converging trees).

\begin{definition}
\label{De2} An out forest $F$ of a digraph $\G$ is called a {\em maximum out forest\/} of $\G$ if $\G$ has no
out forest with a greater number of arcs than in~$F$. An in forest $F$ of a digraph $\G$ is a {\em maximum in
forest\/} of $\G$ if $\G$ has no in forest with a greater number of arcs than in~$F$.
\end{definition}

Obviously, every maximum out forest of $\G$ has the minimum possible number of weak components (out trees);
this number will be called the {\it out forest dimension\/} of the digraph and denoted by~$v$. The number of
arcs in any maximum out forest is obviously $n-v$. The number of weak components of every maximum in forest
will be called the {\it in forest dimension\/} of the digraph and denoted by~$v'$. Obviously, for every
digraph, $v,v'\in\{1\cdc n\}$.

If a digraph $\G_1$ is obtained from $\G$ by the reversal of all arcs, then the out forests in $\G$ naturally
correspond to the in forests in $\G_1$ and vice versa. Therefore, the out forest dimension and in forest
dimension of $\G$ are respectively equal to the in forest dimension and out forest dimension of $\G_1$.

The following proposition states that the dimensions $v$ and $v'$ of a digraph are not connected,{\x} except
for the case where $v=n$ and $v'=n$.

\begin{proposition}
\label{propvv'} $1.$ Let $k,k'\in\{1\cdc n-1\}$. Then there exists a digraph on $n$ vertices such that $v=k$
and $v'=k'$.

$2.$ For every digraph $\G$ on $n$ vertices$,${\x} $v=n\Leftrightarrow v'=n\Leftrightarrow
E(\G)=\emptyset$.{\x}
\end{proposition}

The proofs are given in the Appendix.

Throughout the paper, we mainly deal with diverging forests. However, all the results have counterparts
formulated in terms of converging forests. Simple properties of out forests have been studied
in~\cite{Che.ra1} (Section~3). We do not cite them here and only confine ourselves to the following

\begin{proposition}
\label{prop2} If $i$ and $j$ belong to different trees in a maximum out forest $F$ of a digraph $\G,$ and $j$
is a root in $F,$ then $\G$ contains no paths from $i$ to~$j$.
\end{proposition}

Let us adduce some definitions and results from~\cite{Che.ra1} which are frequently used below.

\setcounter{footnote}{1}
\begin{definition}
\label{De1} A nonempty subset of vertices $K\subseteq V(\G)$ of digraph $\G$ is an {\em undominated
knot\/}\footnote{In~\cite{Fiedler}, undominated knots are called W-bases.} in $\G$ iff all the vertices that
belong to $K$ are mutually reachable and there are no arcs $(w_j,w_i)$ such that $w_j\in V(\G)\setminus K$
and $w_i\in K$.
\end{definition}

Suppose that $\ktil=\cupo^{u}_{i=1}K_i$, where $K_1\cdc K_u$ are all the undominated knots of~$\G$, and
$K_i^{+}$ is the set of all vertices reachable from $K_i$ and unreachable from the other undominated knots.
For any undominated knot $K$ of $\G,$ denote by $\G_K$ the restriction of $\G$ to $K$ and by $\G_{-K}$ the
subgraph with vertex set $V(\G)$ and arc set $E(\G)\setminus E(\G_K)$. For a fixed $K$, $\TT$ will designate
the set of all spanning diverging trees of $\G_K$ and $\PP$ will be the set of all maximum out forests of
$\G_{-K}$. By $\TT^k$, $k\in K$, we denote the subset of $\TT$ consisting of all trees that diverge from $k$,
and by $\PP^{K \rightarrow i}$, $i\in V(\G)$, the set of all maximum out forests of $\G_{-K}$ such that $i$
is reachable from some vertex that belongs to~$K$ in these forests.

By $\FF(\G)=\FF$ and $\FF_k(\G)=\FF_k$ we denote the set of all out forests of $\G$ and the set of all out
forests of $\G$ with $k$ arcs, respectively; $\FF^{i\rightarrow j}_k$ will designate the set of all out
forests with $k$ arcs where $j$ belongs to a tree diverging from~$i$.

\begin{definition}
\label{De3} The matrix $\vj=(\q_{ij})=\si^{-1}Q_{n-v},$ where $\si=\e(\FF_{n-v}),$
$Q_{n-v}=(q\_{ij})=(\e(\FF^{j \to i}_{n-v}))$, will be called the {\it normalized matrix of maximum out
forests\/} of a digraph.
\end{definition}

\begin{theorem}
\label{t1.che.ra1} {Suppose that $\G$ is an arbitrary digraph and $K$ is an undominated knot in~$\G$. Then
the following statements are true.

{\rm 1.}~$\vj$ is a stochastic matrix\/$:${\x} $\q_{ij}\ge0,$ $\;\suml^n_{k=1}\q_{ik}=1,\;$ $i,j=\1n.$

{\rm 2.}~$\q_{ij}\ne 0\;\Leftrightarrow\; (j\in \ktil$ and $i$ is reachable from $j$ in~$\G).$

{\rm 3.}~Suppose that $j\in K.$ For any $i\in V(\G),$
 $\q_{ij}=\e(\TT^j)\e(\PP^{K\rightarrow i})\slash
\e(\FF_{n-v})$. Furthermore$,$ if $i\in K^{+},$ then $\q_{ij}=\q_{jj}=\e(\TT^j)\slash \e(\TT)${\rm.}

{\rm 4.}~$\suml_{j\in K}\q_{jj}=1.$ In particular$,$ if $j$ is an undominated vertex$,$ then $\q_{jj}=1.$

{\rm 5.}~If $j_1,j_2\in K$, then $\q_{\cdot j\_2}=(\e(\TT^{j_2}) \slash\e(\TT^{j_1}))\q_{\cdot j_1},$ i.e.$,$
the $j_1$ and $j_2$ columns of $\q$ are proportional{\rm.}}
\end{theorem}

\begin{theorem}
\label{t2.che.ra1} For every weighted digraph$,$ $\vj$ is idempotent$:$ $\;\q^{2}=\q.$
\end{theorem}

\begin{theorem}
\label{t3.che.ra1} For every weighted digraph$,$ $L\q=\q L=0$.
\end{theorem}

\addtocounter{theorem}{1} {\bf Theorem~4\/} $($a parametric version of the matrix-forest theorem$).$ {\it For
any weighted multidigraph $\G$ with positive weights of arcs and any $\tau>0,$ there exists the matrix
$Q(\tau)=(I+{\tau}L(\G))^{-1}$ and
\begin{equation}
Q(\tau)=\frac{1}{s(\tau)} \suml^{n-v}_{k=0}{\tau}^k Q_k, \label{razlo}
\end{equation}
where
\begin{equation}
\label{raz2} s(\tau)=\suml^{n-v}_{k=0}{\tau}^k\e(\FF_k), \;\; Q_k=(q^k_{ij}), \;\; q^k_{ij}=\e(\FF^{j\to
i}_k), \;\; k=0\cdc n-v, \;\; i,j=\1n.
\end{equation}
}

\begin{definition}
\label{DeQk} The matrix $Q_k,\,${\x} $k=0\cdc n-v,$ will be called the {\it matrix of out forests of $\G$
with $k$ arcs}.
\end{definition}

Theorem~4 represents $(I+\tau L)^{-1}$ via the matrices of out forests with various numbers of arcs.

\begin{definition}
\label{DeQ(tau)} The matrices $Q(\tau)=(I+\tau L)^{-1},\,$ $\tau>0,$ will be called the {\em normalized
matrices of out forests\/} of a digraph.
\end{definition}


In \cite{CheSha97}, the matrices $Q(\tau)=(I+\tau L)^{-1}$ were referred to as the matrices of relative
forest accessibilities of a digraph. In Section~\ref{Sec_Poli}, $Q(\tau)$ are expressed as polynomials of $L$
(Corollary from Theorem~\ref{sumsum7}).

\begin{theorem}
\label{t5.che.le2} For every weighted digraph $\G,$
$\liml_{\tau\to\infty}Q(\tau)=\liml_{\tau\to\infty}(I+\tau\,L)^{-1}=\vj.$
\end{theorem}

\section{Matrices of out forests and transition probabilities of Markov
chains} \label{forest}

It has been shown in~\cite{Che.ra1} that the matrix of Ces\`aro limiting probabilities of a Markov chain{\x}
coincides with the normalized matrix $\vj$ of maximum out forests of any digraph related to this Markov
chain. Now we give a Markov chain interpretation for the normalized matrices of out forests $Q(\tau)$ with
any $\tau>0.$

\bigskip
\addtocounter{definition}{1} {\bf Definition~6\/~\cite{Che.ra1}.} A homogeneous Markov chain with set of
states $\{\1n\}$ and transition probability matrix $P$ is {\it related to a weighted digraph} $\G$ iff there
exists $\alpha\ne0$ such that

\begin{equation}
\label{7.1} P=I-\alpha\,L(\G).
\end{equation}
\smallskip

Let $\G$ be a weighted digraph. Consider an arbitrary Markov chain {\em related to\/}~$\G$ and the following
observation model.
\medskip

{\bf The geometric model of random observation.} {{\x}Suppose that a Bernoulli trial is performed at the
point of time $t=0$ with success probability $q$ $(0<q<1)$. In case of success$,$ $t=0$ becomes the epoch of
observation. Otherwise, Bernoulli trials are performed at $t=1,2,\ldots$---{\x}to the {\x}point of the first
success. This point becomes the epoch of observation.}

\medskip
This model determines a discrete probability distribution $p(k)$ of the epoch of observation on the set
$\{0,1,2,\ldots\}$. This is obviously the geometric distribution (which gives the name of the model) with
parameter~$q$:
\begin{equation}
\label{Geo} p(k)=q(1-q)^k,\quad k=1,2,\ldots
\end{equation}

Consider Markov chain multistep{\x} transitions in a {\em random number of steps}: from the initial state at
$t=0$ to the state at the random epoch of observation distributed geometrically with parameter~$q$.

Suppose that $\Ptil(\aa,q)=\left(\ptil\_{ij}(\aa,q)\right)$ is the matrix of unconditional probabilities for
such multistep transitions: from the initial state to the state at the epoch of observation.

\begin{theorem}
\label{prop1} For any weighted digraph$,$ any $\tau>0$ and any Markov chains related to the weighted
digraph$,$
\[
Q(\tau)=\Ptil(\aa,q)
\]
holds$,$ where
\begin{equation}
\label{qtaual} q=({\tau/\aa}+1)^{-1}.
\end{equation}
\end{theorem}

Theorem~\ref{prop1} provides an interpretation for the normalized{\x} matrix $Q(\tau)$ of out forests{\x} in
terms of Markov chain transition probabilities. Conversely, for any Markov chain, the transition
probabilities in the geometric observation model can be interpreted in terms of diverging forests of the
corresponding digraphs.

The following corollary stresses the arbitrariness of Markov chains in Theorem~\ref{prop1}.

\medskip
\begin{coroll}
{\bf 1 from Theorem~\ref{prop1}} {For every Markov chain$,$ every success probability $q\in\,]0,1[$ in the
geometric observation model$,${\x} and every digraph related to the Markov chain$,$
\[
\Ptil(\aa,q)=Q(\tau)
\]
holds$,$ where $\tau=(q^{-1}-1)\aa.$ }
\end{coroll}
\medskip

\begin{coroll}
{\bf 2 from Theorem~\ref{prop1}}{
\begin{equation}
\label{Cor1} \lim_{q\to+0}\Ptil(\aa,q)=\vj=\lim_{k\to\infty} \frac{1}{k} \suml_{p=0}^{k-1} P^{k}.
\end{equation}
}
\end{coroll}

By Corollary~2 from Theorem~\ref{prop1}, at a vanishingly small success probability $q$, the transition
probabilities in the geometric observation model are given{\x} by the matrix $\vj$ of maximum out forests of
any weighted digraph to which this chain is related.

\section{Representations of forest matrices via the Kirchhoff matrix
and their consequences} \label{Sec_Poli}

In this section, we represent the matrices $Q_k$ of out forests with $k$ arcs as polynomials of the Kirchhoff
matrix $L$ (Theorem~\ref{teo.allk}). This allows one to obtain{\x} alternative proofs of Theorems~2--4 and to
represent the matrix $Q(\tau)=(I+\tau L)^{-1}$ as a polynomial of $L$ (Theorem~\ref{sumsum7}).
Proposition~\ref{recur} gives an easy way to calculate $Q_k,\,$ $k=1\cdc n-v,${\x} and~$\q$.

By $\si\_{k}$ we denote the total weight of all out forests of $\G$ with $k$ arcs: $\si\_{k}=\e(\FF_{k}),\,$
$k=0\cdc n-v.$

\begin{proposition}
\label{pro.allk} For any weighted digraph and any $k=0\cdc n-v,$
\begin{equation}
\label{QQQ} Q_{k+1}=\si\_{k+1}\!I-L Q_{k}.
\end{equation}
\end{proposition}

Observe that since the weight of the empty set is 0, we have $Q_{n-v+1}=0$ and $\si\_{n-v+1}=0$.

Taking the traces on the left-hand side and the right-hand side of (\ref{QQQ}) and using the fact that
\begin{equation}
\label{TrQ} \tr{\x}(Q_{k})=(n-k)\si\_{k},\quad k=0\cdc n-v+1
\end{equation}
(because every out forest with $k$ arcs has $n-k$ roots), we deduce
\begin{equation}
\label{si_k+1} \si_{k+1}=\frac{\tr{\x}(LQ_k)}{k+1},\quad k=0\cdc n-v.
\end{equation}

Substituting (\ref{si_k+1}) in (\ref{QQQ}) provides

\begin{proposition}
\label{recur} For every weighted digraph$,${\x}
\begin{equation}
\label{QQ+} Q_{k+1}=\frac{\tr{\x}(LQ_k)}{k+1}I-LQ_{k},\quad k=0\cdc n-v.
\end{equation}
\end{proposition}

Identity (\ref{QQ+}) enables one to recursively determine the matrices $Q_k,\,$ $k=0\cdc n-v,${\x} and $\q$,
starting with $Q_0=I.$ Note that this procedure essentially coincides with Faddeev's algorithm \cite{Faddeev}
for the computation of the characteristic polynomial as applied to~$L$. Thus, the matrices involved in
Faddeev's method are precisely $Q_k$.

From Proposition~\ref{pro.allk}, it follows
\begin{theorem}
\label{teo.allk} For any weighted digraph and any $k=0\cdc n-v,$
\begin{equation}
\label{psibek} Q_{k}=\suml_{i=0}^{k} \si\_{k-i} {(-L)}^{i}.
\end{equation}
\end{theorem}

\begin{coroll}
{\bf 1 from Theorem~\ref{teo.allk}} {For every weighted digraph$,$ matrices $Q_k,$ $k=0\cdc n-v,$ commute
with all matrices with which $L$ commutes$,$ in particular$,$ with $L,$ $\vj,$ $Q(\tau),$ and each other. }
\end{coroll}

\begin{lemma}
\label{sumstr} For any $k=0\cdc n-v,$ every row sum of $L Q_{k}$ is~$0$.
\end{lemma}

From Proposition~\ref{pro.allk} and Lemma~\ref{sumstr}, it follows

\begin{proposition}
\label{allmatrbek} The matrices $L Q_{k},$ $k=0\cdc n-v,$ are the Kirchhoff matrices of some weighted
digraphs.
\end{proposition}

\begin{coroll}
{\bf 2 from Theorem~\ref{teo.allk}} {For any weighted digraph$,$ $LQ_{n-v}=Q_{n-v}L=0.$}
\end{coroll}
\smallskip

In view of Definition~\ref{De3}, this corollary is equivalent to Theorem~\ref{t3.che.ra1}. Thus,{\x} we get a
new proof of this theorem.

Consider the matrices
\begin{equation}
\label{*} \q_k=\si_{k}^{-1}Q_k,\quad k=0{\x}\cdc n-v.
\end{equation}
In particular, $\q_0=I$ and $\q_{n-v}=\q$.

Making use of the last corollary, we obtain

\medskip
\begin{coroll}
{\bf 3 from Theorem~\ref{teo.allk}} {For any $k\in\{\1n-v\},\,$ $\q_k\q=\q\q_k=\q.$ In particular$,$
$\q_{n-v}\vj=\vj^2=\vj.$ Moreover$,$ $Q(\tau)\q=\q Q(\tau)=\q$ for every $\tau>0.$ }
\end{coroll}
\smallskip{\x}

This corollary provides a new proof of Theorem~\ref{t2.che.ra1}.
\smallskip

By virtue of Proposition~\ref{pro.allk} and Corollary~1 from Theorem~\ref{teo.allk}, the matrices $\q_{k}$
are connected{\x} as follows:
\begin{equation}
\label{j_k+1} \q_{k+1}=I-\frac{\si\_k}{\si\_{k+1}}\q_{k} L, \quad k=0\cdc n-v-1,
\end{equation}
and,{\x} by Lemma~\ref{sumstr}, each their row is unity.{\x} The entries of $\q_k$ are nonnegative by
definition, thus, we obtain

\begin{proposition}
\label{propQk1} For every weighted digraph $\G,$ matrices $\q_k,\,$ $k=0\cdc n-v,${\x} are stochastic.
\end{proposition}

Completing Proposition~\ref{pro.allk} with the obvious equality $Q_{0}=I=\si\_0 I$ gives
\begin{equation}
\label{system1} \cases{
   Q_{0}=\si\_0 I,                       \cr
   Q_{1}+L Q_{0}=\si\_1 I,              \cr
   \ldots\ldots\ldots\ldots\ldots\ldots  \cr
   Q_{n-v}+L Q_{n-v-1}=\si\_{n-v} I.     \cr
}
\end{equation}

Add up these equations{\x} and, using Corollary~2 from Theorem~\ref{teo.allk}, substitute $(I+L)Q_{n-v}$ for
$Q_{n-v}$:
\[
(I+L)Q_{0}+(I+L)Q_{1}+\ldots+(I+L)Q_{n-v} =\Bigl(\suml_{k=0}^{n-v}\si\_k\Bigr)I.
\]

Making use of{\x} the nonsingularity of $I+L$ (Theorem~4) and the notation
$s=\e(\FF)=\suml_{k=0}^{n-v}\si\_k$, we obtain
\begin{equation}
\label{iden1} \suml_{k=0}^{n-v}Q_k = s(I+L)^{-1},
\end{equation}
which provides

\medskip
\begin{coroll}
{\bf from Proposition~\ref{pro.allk}} {For any weighted digraph$,${\x}
\[
Q(1)=(I+L)^{-1}=s^{-1}\suml_{k=0}^{n-v}Q_k.
\]
}
\end{coroll}
\smallskip

This statement coincides with the matrix-forest theorem for digraphs~\cite{CheSha97} and with Theorem~4 in
the case of $\tau=1$. Accordingly, we obtain a new proof of the matrix-forest theorem.

By means of Theorem~\ref{teo.allk}, the matrices $Q(\tau)={(I+\tau L)}^{-1}=s^{-1}\suml_{k=0}^{n-v}\tau^k
Q_k${\x} including $Q(1)={(I+L)}^{-1}${\x} can be represented as polynomials of~$L$.

\begin{theorem}
\label{sumsum7} For any weighted digraph$,${\x}
\begin{equation}
\label{eqth7} Q(1)={(I+L)}^{-1}=s^{-1}\suml_{i=0}^{n-v}s\_{n-v-i}(-L)^i,
\end{equation}
where $s\_k=\suml_{j=0}^k\si\_j$ is the total weight of out forests of $\G$ with at most $k$ arcs$,$ $k=0\cdc
n-v$.{\x}
\end{theorem}

\begin{coroll}
{\bf from Theorem~\ref{sumsum7}} {For any weighted digraph and any $\tau>0,$
\begin{equation}
\label{eqco7} Q(\tau)={(I+\tau L)}^{-1}= s^{-1}(\tau)\suml_{i=0}^{n-v}s\_{n-v-i}(\tau)(-\tau L)^i,
\end{equation}
where $s\_k(\tau)=\suml_{j=0}^{k}\tau^j\si\_j,\:$ $k=0\cdc n-v$.{\x} }
\end{coroll}
\medskip

Note that $s(I+L)^{-1}$ is the adjugate (the transposed matrix of cofactors) of $I+L${\x}; $s(\tau)\,(I+\tau
L)^{-1}$ is the same for $I+\tau L${\x}. Theorem~\ref{sumsum7} and the above corollary provide
representations for these matrices as polynomials of~$L$:
\begin{eqnarray}
\label{iden2} s\,(I+L)^{-1} &=&\suml_{i=0}^{n-v}s\_{n-v-i}(-L)^i,\cr \label{iden3} s(\tau)\,(I+\tau L)^{-1}
&=&\suml_{i=0}^{n-v}s\_{n-v-i}(\tau)\,(-\tau L)^i.
\end{eqnarray}

\begin{remark}
Since $L Q_k$ the is Kirchhoff matrix of some weighted digraph (Proposition~\ref{allmatrbek}), all its
principal minors are nonnegative (by Theorem~6 in~\cite{Fiedler}). Therefore, all $L Q_k$ are singular
$M$-matrices (see, e.g., item~(A1) of Theorem~4.6 in~\cite{Abraham.Robert}). Alternatively, this can be
concluded from the nonnegativity of the real parts of the eigenvalues (see Proposition~\ref{Gosha} below) and
the nonpositivity of off-diagonal elements of $L$ (item~(F12) of Theorem~4.6 in \cite{Abraham.Robert}). It
follows from the representation $\si\_{k+1}I-Q_{k+1}=LQ_k$ (Proposition~\ref{pro.allk}) of the singular
$M$-matrix $L Q_k$ that $\si\_{k+1}=\rho(Q_{k+1})$, i.e., $\si\_{k+1}$ is the spectral radius of $Q_{k+1}$,
$k=0\cdc n-v-1${\x}. This also follows from Proposition~\ref{propQk1} (see~(\ref{*})).
\end{remark}

\section{On some linear transformations related to digraphs}
\label{line}

For a matrix $A\in\R^{n\times n}$, by ${\bf A}$ we denote the linear transformation ${\bf A}:\R^n\to\R^n$
induced by $A$ with respect to the standard basis of $\R^n$: ${\bf A}({\bf x})=A{\bf x}$. $\RR({\bf A})$ and
$\NN({\bf A})$ will designate the range and the null space of ${\bf A}$, respectively.{\x}

As has been seen in~\cite{Che.ra1}, the dimensions of $\RJ$ and $\RLT$ are $v$ and $n-v$, respectively.
Furthermore, $\RLT\cap\RJ=\{{\bf 0}\}$ and, since the dimensions of $\RLT$ and $\RJ$ sum to $n$, $\R^n$
decomposes to the direct sum of $\RLT$ and $\RJ$:{\x}
\begin{equation}
\label{pr.summa} \R^n=\RLT\dot+\RJ.{\x}
\end{equation}

Since $L\!\vj=0$ (Theorem~\ref{t3.che.ra1}), we get $\NL=\RJ$ and $\NJT=\RLT$, thus, the sum (\ref{pr.summa})
is orthogonal.{\x}

Similarly, in view of $\vj\! L=0$, the orthogonal decomposition
$$
\R^n=\RL\dot+\RJT{\x}
$$
holds along with $\RL\cap\RJT=\{{\bf 0}\}$, $\NJ=\RL$, and $\NLT=\RJT$.{\x}

In accordance with~(\ref{pr.summa}), every vector ${\bf u}\in\R^n$ is uniquely represented as ${\bf u}={\bf
u}\_1+{\bf u}\_2,$ where ${\bf u}\_1\in \RLT=\NJT$ and{\x} ${\bf u}\_2\in\RJ=\NL$. For every ${\bf u}\ne{\bf
0}$, we have ${(L+\vj}^{\intercal}) {\bf u}={(L+\vj}^{\intercal}){\bf u}\_1+$ ${(L+\vj}^{\intercal}) {\bf
u}\_2$ $=L{\bf u}\_1+{\vj}^{\intercal}{\bf u}\_2.$ If $L{\bf u}\_1+{\vj}^{\intercal}{\bf u}\_2=0$, then,
since $\RL\cap\RJT=\{{\bf 0}\}$, we have $L{\bf u}\_1={\vj}^{\intercal}{\bf u}\_2={\bf 0}$, whence, by
$\NL\cap\NJT=\{{\bf 0}\}$, ${\bf u}\_1={\bf u}\_2={\bf 0}$ results. Therefore, the dimension of the range
(rank) of ${\bf Z=L+\vj}^{\intercal}$ is~$n$. Thus, we obtain{\x}

\begin{theorem}
\label{invert.Ltj} For any weighted digraph $\G,$ the matrix $Z=L+\vj^{\intercal}$ is nonsingular.
\end{theorem}

{\x}{\x}

We will also need the nonsingularity of $L+\vj$.

\begin{theorem}
\label{fulllj} For any weighted digraph $\G,$ the matrix $L+\q$ is nonsingular.
\end{theorem}

\begin{coroll}
{\bf from Theorem~\ref{fulllj}} {For any weighted digraph and any $\aa\ne 0,$ the matrix $L+\aa\vj$ is
nonsingular.}
\end{coroll}
\smallskip

It follows from $\vj^2=\vj$ (Theorem~\ref{t2.che.ra1}) that every{\x} nonzero columns of $\vj$ is an
eigenvector of $\vj$ associated with the eigenvalue~1. Hence, for any ${\bf u}\in\RJ$,{\x} ${\bf\vj}({\bf
u})={\bf u}$ holds, therefore, $\RJ$ is exactly the subspace of fixed vectors of ${\bf\vj}$.{\x}

\section{The Moore-Penrose and group inverses of the Kirchhoff matrix}
\label{pseudo2}

In this section, we obtain some expressions for the Moore-Penrose generalized inverse and the group inverse
of the Kirchhoff matrix~$L$. The {\it Moore-Penrose generalized inverse\/} of a rectangular complex matrix
$A$ is the unique matrix $X$ such that

(1)~$AXA=A,$

(2)~$XAX=X,$

(3)~$(AX)^*=AX,$ 

(4)~$(XA)^*=XA,$ 

\noindent where $(AX)^*$ and $(XA)^*$ are the conjugate transposes (Hermitian adjoints) of $AX$ and $XA$,
respectively{\x}.

For any matrix $A,$ the Moore-Penrose generalized inverse, $A^+$, does exist{\x} and is unique. If $A$ is
nonsingular, then $A^+$ coincides with~$A^{-1}$.

The Moore-Penrose inverses are of theoretical and practical interest. The latter is because $A^+$ provides
the normal pseudosolution of the inconsistent equation $A{\bf x}={\bf b}$: it is ${\bf x}=A^+{\bf b}$. The
normal pseudosolution is a vector of the minimum length that minimizes the length of $A{\bf x}-{\bf b}$ (the
minimum norm least-squares solution). As applied to Laplacian matrices, such solutions, among others, were
considered for some preference aggregation problems (more specifically, estimation from paired comparisons)
\cite{CheSha99}, in constructing geometrical representations for systems modelled by graphs \cite{Hall}, in
the analysis of social networks, and cluster analysis.

The group inverses are no less important (see, e.g.,~\cite{CampMey}). A~matrix $X$ is the {\em group
inverse\/} of a square matrix $A$, if $X$ satisfies the conditions (1) and (2) in the definition of
Moore-Penrose generalized inverse and also

(5) $AX=XA$.

The group inverse of $A$ is denoted by~$\Apr$. Generally, group inverses need not exist, but if such a matrix
exists, then it is unique, but $A^+=\Apr$ is not necessary.

If $L$ is symmetric (in particular, this is the case for symmetric digraphs, which can be identified with
undirected graphs{\x}), then the matrix $(L+\aa\vj)^{-1}-\aa^{-1}\vj$ (with any $\aa>0$) is~\cite{CheSha981}
the Moore-Penrose generalized inverse and the group inverse of~$L$. Moreover, the latter is true for every
digraph.

\begin{theorem}
\label{groupinv} For every weighted digraph and any $\aa\ne0,$
\begin{equation}
\label{qtil} \Lpr=(L+\aa\vj)^{-1}-\aa^{-1}\vj
\end{equation}
and
\[
\Lpr L=L\Lpr=I-\q.
\]
\end{theorem}

As well as in the case of undirected graphs, $\Lpr=(\lpr_{ij})$ can be obtained via a passage to the limit.

\begin{proposition}
\label{qtiLim} For every weighted digraph$,${\x}
\[
\Lpr=\liml_{\tau\to\infty}\tau\Big(Q(\tau)-\vj\Big).
\]
\end{proposition}

We now express $\Lpr$ in terms of the normalized matrices of out forests $\q_{n-v-1}$ and $\q_{n-v}=\q$
(see~(\ref{*})). The following proposition is an analogue of Theorem~3 in~\cite{CheSha981}.{\x}

\begin{proposition}
\label{topol} For every weighted digraph$,$
\[
\Lpr=\frac{\si_{n-v-1}}{\si_{n-v}}\left(\q_{n-v-1}-\q\right).
\]
\end{proposition}

Because of the nonsymmetry of $I-\q=\Lpr L=L\Lpr$, $\Lpr$ is not generally the Moore-Penrose generalized
inverse of $L$ for digraphs. To obtain an explicit formula for $L^+,$ consider the matrix
$Z=L+\vj^{\intercal}${\x} which, by Theorem~\ref{invert.Ltj}, is nonsingular. Using the identity $L\vj=0$
(Theorem~\ref{t3.che.ra1}), we obtain
\[
{(Z^{\intercal})}^{-1}Z^{-1} =(ZZ^{\intercal})^{-1} =(\vj^{\intercal}\!\vj+LL^{\intercal})^{-1}.
\]

\begin{lemma}
\label{commut} For every weighted digraph$,$ ${(ZZ^{\intercal})}^{-1}$ commutes with $LL^{\intercal}$ and
$\vj^{\intercal}\!\vj$.
\end{lemma}

Matrices $L L^{\intercal}$, $\vj^{\intercal}\!\vj$, and ${(ZZ^{\intercal})}^{-1}$ are symmetric. The product
of two symmetric matrices is symmetric iff they are commuting~\cite{HoJo}. This implies the following
corollary.

\medskip
\begin{coroll}
{\bf from Lemma~\ref{commut}} {For every weighted digraph$,$ the matrices
$LL^{\intercal}{(ZZ^{\intercal})}^{-1}$ and
$\vj^{\intercal}\!\vj{(ZZ^{\intercal})}^{-1}$ are symmetric. }
\end{coroll}
\smallskip{\x}

These facts are useful for the proof of the following theorem.

\begin{theorem}
\label{pseudoor} For every weighted digraph$,$ the matrix $L^{\intercal}{(ZZ^{\intercal})}^{-1}=L^{\intercal}
(\vj^{\intercal}\!\vj+L L^{\intercal})^{-1}$ is the Moore-Penrose generalized inverse of~$L$.
\end{theorem}

\section{On the Ger\v{s}gorin region and annihilating polynomials
for the Kirchhoff matrix}

By the Ger\v{s}gorin theorem (see, e.g.,~\cite{HoJo}), the eigenvalues of a matrix $A$ belong to the union
$G(A)$ of $n$ discs:
\begin{equation}
\label{Gershgorn} G(A)=\cupo_{i=1}^{n}\Bigl\{z\in\C\:\Big\vert\;|z-a_{ii}|\le R'_i(A)\Bigr\},
\end{equation}
where $\C$ is the complex field and $R'_i(A)=\suml_{j \ne i}|a_{ij}|,$ $i=\1n,$ are the deleted absolute row
sums of~$A$.

Since $R'_i (L)=\l_{ii}$ holds, (\ref{Gershgorn}) can be represented as follows:
\begin{equation}
\label{GershgornL} G(L)=\cupo_{i=1}^{n}\Bigl\{z\in\C\;\Big\vert\:|z-\l_{ii}| \le\l_{ii}\Bigr\}.
\end{equation}

Hence, we have

\begin{proposition}
\label{Gosha} {\rm (1)} The real part of every eigenvalue of $L$ is nonnegative$:$ every Ger\v{s}gorin disc
belongs to the right coordinate half-plane$;$

{\rm (2)} the intersection of all Ger\v{s}gorin discs contains zero$;$

{\rm (3)} $G(L)=\Bigl\{(z+1) \maxl_{1\le i\le n} \l_{ii}\:\Big\vert\; |z| \le 1\Bigr\}$.
\end{proposition}

Obviously, the intersection of all Ger\v{s}gorin discs consists of zero {\x}iff the digraph contain an
undominated vertex.

Consider the characteristic polynomial of $L$:
\[
p\_L(\la)=\suml_{i=0}^{n}{(-1)}^i E_{i}(L)\la^{n-i},
\]
where $E_{i}(L)$ is the sum of all principal minors of order~$i$. By Theorem~6 in \cite{Fiedler},
$E_{i}(L)=\si\_i$ for every $i=\1n$. Since every principal minor of order greater than $n-v$ is zero, we have
 $p\_L(\la)={(-1)}^{n-v}\la^v \suml_{i=0}^{n-v}
\si\_{i}{(-\la)}^{n-v-i}=\la^v\suml_{i=0}^{n-v}(-1)^{i} \si\_{i}{\la}^{n-v-i}$.

\begin{proposition}
\label{annul} $p'_L(\la)=\la \suml_{i=0}^{n-v}\si\_{n-v-i}{(-\la)}^{i}$ is an annihilating polynomial for
$L$.
\end{proposition}

\section{Accessibility via forests and dense forests in digraphs}
\label{dostup}

\subsection{Forest accessibility}
\label{dostup1} The entries of $Q(\tau)$ measure the proximity of the vertices of an undirected
multigraph~\cite{CheSha97,CheSha981}. The matrix $\vj^{\intercal}=\lim_{\tau\to \infty}Q^{\intercal}(\tau)$
was analyzed in \cite{Che.ra1} as the matrix of limiting accessibilities of a multidigraph. Here, we study
the matrix $P\_1(\tau)=Q^{\intercal}(\tau)$ with $\tau>0$ as an accessibility measure for digraph vertices.
By Theorem~4, the $(i,j)$-entry of this matrix is the total weight of out forests that ``connect'' $i$ with
$j$ in the digraph where the weights of all arcs are multiplied by~$\tau$. Along with $P\_1(\tau)$,{\x} we
consider the matrix of in forests $P\_2(\tau)$. Its $(i,j)$-entry is the total weight of in forests (of{\x}
the modified digraph) where $j$ is a sink and $i$~belongs to a tree converging to~$j$.

The following definition is formulated for an arbitrary vertex accessibility measure (formally, every square
matrix of order $n$ or, more precisely, the corresponding matrix-valued function of a digraph can be
considered as such a measure). A measure $P\_2$ is called to be {\em dual\/} to a measure $P\_1$ if under the
reversal of all arcs in an arbitrary digraph (provided that the weights of the arcs are preserved),{\x} the
matrix of $P\_2$ for the modified digraph coincides with $P_1^{\intercal}$ calculated for the initial
digraph. It follows from this definition that $P\_2$ is dual to $P\_1$ if and only if $P\_1$ is dual to
$P\_2.$ In~\cite{CheSha981}, three self-dual{\x} accessibility measures were studied.

Let us check the satisfaction of the characteristic conditions listed below for $P\_1(\tau)$ and
$P\_2(\tau)$. Triangle inequality for accessibility measures requires the symmetry of the corresponding
matrix (see, e.g.,~\cite{CheSha982}). For that reason, we will check this condition for
$P\_3(\tau)=(P\_{1}(\tau)+P\_{2}(\tau)+ P_{1}^{\intercal}(\tau)+P_{2}^{\intercal}(\tau))/4$.

\axiom{Nonnegativity} {For any digraph $\G,$ $\;p_{ij}\ge 0,\;\:i,j\in V(\G)$.}

\axiomt{Diagonal maximality} {For any digraph $\G$ and any distinct $i,j\in V(\G),$

{\rm (1)} $p\_{ii}>p\_{ij}$ and

{\rm (2)} $p\_{ii}>p\_{ji}$ hold.}

\axiomt{Disconnection condition} {For any digraph $\G$ and any $i,j\in V(\G),\;$ $p\_{ij}=0$ if and only if
$j$ is unreachable from~$i$.}

\axiomt{Triangle inequality for accessibility measures} {For any digraph $\G$ and any $i,j,k\in V(\G),$
$p\_{ij}+p\_{ik}-p\_{jk}$ $\le p\_{ii}$ holds. If, in addition, $j=k$ and $i\ne j$, then the inequality is
strict.}

\axiomt{Transit property} {For any digraph $\G$ and any $i,k,t\in V(\G),$ if $\G$ includes a path from $i$ to
$k,$ $i\ne k\ne t,$ and every path from $i$ to $t$ contains $k,$ then {\rm (1)} $p\_{ik}>p\_{it};$ {\rm (2)}
$p\_{kt}>p\_{it}$.}

\axiomt{Monotonicity} {Suppose that the weight of some arc $\e_{kt}^p$ in a digraph $\G$ increases.
Then{\rm:}
\parr
{\rm (1)} $\D p\_{kt}>0${\x} and for any $i,j\in V(\G),$ $(i,j)\ne (k,t)$ implies $\D p\_{kt}>\D p\_{ij};$
\parr
{\rm (2)} For any $i\in V(\G),$ if there is a path from $k$ to $t,$ and each path from $k$ to $i$ includes
$t,$ then $(a)$ $\D p\_{kt}>\D p\_{ki}$ and $(b)$ $\D p\_{ki}>\D p\_{ti};$

{\rm (3)} For any $i\in V(\G),$ if there is a path from $i$ to $k$ and every path from $i$ to $t$ includes
$k,$ then $(a)$ $\D p\_{kt}>\D p\_{it}$ and $(b)$ $\D p\_{it}>\D p\_{ik}.$}

The results of testing $P\_1(\tau),$ $P\_2(\tau),$ and $P\_3(\tau)$ are collected{\x} in the following
proposition.

\begin{proposition}
\label{otledostup} The measures $P\_1(\tau)$ and $P\_2(\tau)$ are {\x}dual to each other for every $\tau>0$.
They satisfy nonnegativity$,$ reversal property$,$ disconnection condition$,$ the first part of item~$1,$ and
item~$2$ of monotonicity. Moreover$,$ $P\_1(\tau)$ satisfies items~{\rm 1} of diagonal maximality and transit
property$;$ $P\_2(\tau)$ satisfies items~{\rm 2} of these conditions. With respect to the remaining
statements of monotonicity$,$ $P\_1(\tau)$ satisfies items $2$ and $3b$, whereas $P\_2(\tau)$ satisfies items
$3$ and~$2b,$ and they both violate the second part of item~$1.$ Furthermore$,$ $P\_1(\tau)$ breaks item $3a$
of monotonicity and item~$2$ of transit property$,${\x} whereas $P\_2(\tau)$ breaks item~$2a$ of monotonicity
and item~$1$ of transit property. Triangle inequality for $P\_3(\tau)$ is not satisfied.

\end{proposition}

As was noted in~\cite{Che.ra1}, the limiting accessibility $P=\vj^{\intercal}$ of a digraph does not
completely correspond to the general concept of proximity. Notice that {\it disconnection condition}, which
is satisfied for the limiting accessibility in one side only, is completely fulfilled for $P\_1(\tau)$ and
$P\_2(\tau)$. Moreover, $P\_1(\tau)$ and $P\_2(\tau)$ obey a number of conditions which are satisfied by the
limiting accessibility in the nonstrict form only.\footnote{By the nonstrict form of a condition we mean the
result of substituting nonstrict inequalities ($\ge$ and $\le$) for the strict ones ($>$ and $<$) in it.}

\subsection{Accessibility via dense forests}
\label{dostup2} Now we consider a measure which is intermediate between the limiting accessibility (which{\x}
depends on $Q_{n-v}$ only) and the forest accessibility $Q(\tau)$ (which is a weighted sum of all
matrices~$Q_k$). This new measure is determined by the matrices $Q_{n-v-1}$ and $Q_{n-v}$ (or, equivalently,
by the matrices $\q_{n-v-1}$ and $\q_{n-v}=\q$, which also determine $\Lpr$ as stated in
Proposition~\ref{topol}). This measure can be also obtained by the inversion of $L+\aa\vj$ with some values
of~$\aa$.

Thus, consider the matrices $R(\aa)=(r_{ij})=(L+\aa\vj)^{-1}$ with $\aa>0$. Using Theorem~\ref{groupinv} and
Proposition~\ref{topol}, we have
\begin{equation}
\label{f21} (L+\aa\vj)^{-1}=\Lpr+\aa^{-1}\vj =\frac{\si_{n-v-1}}{\si_{n-v}}\q_{n-v-1}
+\left(\aa^{-1}-\frac{\si_{n-v-1}}{\si_{n-v}}\right)\q.
\end{equation}

If $0<\aa<\frac{\si_{n-v}}{\si_{n-v-1}}$, then, by (\ref{f21}), $(L+\aa\oj)^{-1}$ is the sum of $Q_{n-v-1}$
and $Q_{n-v}$ with positive coefficients. Spanning rooted forests (of an undirected multigraph) with $n-v$ or
$n-v-1$ arcs are called in \cite{CheSha981} {\em dense forests}, and the undirected counterpart of the
accessibility measure{\x} (\ref{f21}) with $0<\aa<\frac{\si_{n-v}}{\si_{n-v-1}}$ is called {\em accessibility
via dense forests}.

Consider two accessibility measures for digraphs: $P\_1(\aa)=R^{\intercal}(\aa)$, {\em accessibility via
dense diverging forests\/} and $P\_2(\aa)$, {\em accessibility via dense converging forests}.

An important property of the set of dense diverging forests is as follows.

\begin{proposition}
\label{p1.2} $1.$ For any vertex $i\in V(\G),$ there exists an out forest in $\FF_{n-v-1}$ where $i$ is a
root. $2.$ For any path $($chain subgraph\/{\x}$)$ in $\G,$ there exists an out forest in
$\FF_{n-v-1}\cup\FF_{n-v}$ that contains this path.
\end{proposition}

A similar proposition is true for converging forests. At the same time, the set of maximum out forests
$\FF_{n-v}$ and the set of maximum in forests do not have this property. For example, on Fig.~1 in
\cite{Che.ra1}, no maximum out forest contains arc~$(4,2)$.

We now{\x} test $P\_1(\aa)$ and $P\_2(\aa).$ Similar to the previous consideration, triangle inequality for
accessibility measures will be checked for the index
$P\_3(\aa)=(P\_1(\aa)+P_1^{\intercal}(\aa)+P\_2(\aa)+P_2^{\intercal}(\aa))/4$, since this inequality requires
the symmetry of the corresponding matrix.

\begin{proposition}
\label{pogushe} For any $\aa\in\:]\,0,\,\si_{n-v}/\si_{n-v-1}\/[,$ the measures $P\_1(\aa)$ and $P\_2(\aa)$
are dual to each other. They satisfy nonnegativity and disconnection condition. Moreover$,$ the{\x} nonstrict
versions of items~{\rm 1} of diagonal maximality and transit property are satisfied by $P\_1(\aa),$ and
items~{\rm 2} of these conditions by $P\_2(\aa)$. Both measures violate monotonicity. Triangle inequality for
accessibility measures is not true for $P\_3(\aa)$.
\end{proposition}

\section*{Conclusion}

The normalized matrices of out forests are stochastic and determine the transition probabilities in the
geometric observation model applied to the Markov chains related to the digraph under consideration. Various
matrices of forests can be represented by simple polynomials of the Kirchhoff matrix. The Moore-Penrose
generalized inverse $L^+$ and the group inverse $\Lpr$ of the Kirchhoff matrix $L$ can be explicitly
represented via $L$ and the normalized matrix $\q$ of digraph's maximum out forests. The matrices of
diverging and converging forests characterize the pairwise accessibility of vertices. These and other results
enable one to consider the matrices of spanning forests as a useful tool for the analysis of digraph's
structure.

\section*{\sl\hfill Appendix}

\propro{\ref{propvv'}}{ 1. At first, let{\x} $k\le k'<n.$ $\G$ is{\x} constructed as follows. We draw a
diverging star rooted at the first vertex and having $k'-k$ leaf vertices. Also, we draw a path diverging
from the root of the star and containing, in addition to the root, $n-k'\ge1$ vertices that are not included
in the star. The remaining $k-1$ vertices are left isolated. Then $v=1+(k-1)=k$ and $v'=1+(k'-k)+(k-1)=k'$,
as required.

If $k'\le k<n,$ then we draw a star converging to the first vertex and having $k-k'$ other vertices. Also, we
draw a path diverging to the center of the star and containing, in addition to the sink, $n-k\ge1$ vertices
that are not included in the star. The remaining $k'-1$ vertices are left isolated. As well as in the first
case, we have $v=k$ and $v'=k'$. The second statement is obvious.}

\prothe{\ref{prop1}}{ Since the spectral radius of $P$ is 1, we have
\[
\suml_{k=0}^{\infty}\Big((1-q)P\Big)^k=\Bigl(I-(1-q)P\Bigr)^{-1}.
\]

Using the formula of total probability and equations (\ref{7.1})--%
(\ref{qtaual}), we obtain
\begin{eqnarray}
\Ptil(\aa,q) &\!=&\!\suml_{k=0}^{\infty}p(k)P^k=\!\suml_{k=0}^{\infty}
       q(1-q)^kP^k=q\Big(I-(1-q)P\Big)^{-1}\\
&\!=&\!q\Big(I-(1-q)(I-\aa L)\Big)^{-1}=q\Big(qI+(1-q)\aa
       L\Big)^{-1}\nonumber\\
&\!=&\!\left(I+\frac{(1-q)\aa}{q}L\right)^{-1}=(I+\tau
       L)^{-1}=Q(\tau).
\nonumber
\end{eqnarray}
}

\propro{\ref{pro.allk}}{
In \cite{Che.ra1}, we used the notion of {\x}weighted 2-digraph: it is a multidigraph with arc multiplicities
no more than two. The weight of a 2-digraph is the product of the weights of its arcs. For a weighted digraph
$\h$ and its vertices $u,w\in V(\h)$, by $\,\h+(u,w)$ we denote the 2-digraph with vertex set $V(\h)$ and the
arc multiset obtained from $E(\h)$ by the increment of the multiplicity of $(u,w)$ by~1. Similarly, if $\h$
is a 2-digraph and $u,w\in V(\h)$, then by $\h'=\h-(u,w)$ we denote the 2-digraph that differs from $\h$ in
the multiplicity of $(u,w)$ only: $n'((u,w))=\max(n((u,w))-1,0).$

Now introduce the following notation. Let $\FF_{k}^{j\to s}+(\l,i)= \{F_{k}^{j\to s}+(\l,i)\mid\;F_{k}^{j\to
s}\in\FF_{k}^{j\to s}\}$. For all $i\ne j,$ by $\FF_{k,i}^{j\to s}$ denote the set of out forests with $k$
arcs where $i$ is a root and $s$ belongs to a tree diverging from~$j$. Obviously, in such digraphs, $s$ is
unreachable from $i$ whenever $i\ne j$. By $\FF_{k,{\bar i}}^{j\to s}$ we denote the set of all forests with
$k$ arcs where $i$ {\it is not\/} a root, and $s$ belongs to a tree diverging from~$j$. These definitions
induce{\x} the matrices $(q^{k}_{sj,{\bar i}})$ and $(q^{k}_{sj,i})$ with elements
\begin{eqnarray*}
q^{k}_{sj,{\bar i}}
&=&\e(\FF_{k,{\bar i}}^{j\to s}),\\
q^{k}_{sj,i} &=&\e(\FF_{k,i}^{j\to s}), \;\; s,j=\1n.
\end{eqnarray*}

Denote by $q^{k}_{{\bar i}j}$ the total weight of all diverging forests where $i$ is unreachable from~$j$.
For all $i,j\in V(\G),$ $q^{k}_{jj}=q^{k}_{ij}+q^{k}_{{\bar i}j}$ holds.

\begin{lemma}
\label{matr1} For any weighted digraph and all $i,j=\1n,$ we have{\rm:}

{\rm (1)} If $s\ne p${\rm,} $i\ne j,$ and $(s,i),(p,i)\in E(\G)${\rm,} then $(\FF_{k,i}^{j\to s}+(s,i))\capo
\,(\FF_{k,i}^{j\to p}+(p,i))=\emptyset;$

{\rm (2)} $\cupo_{(s,i)\in E(\G)} (\FF_{k,i}^{j\to s}+(s,i))= \FF_{k+1}^{j\to i}.$
\end{lemma}

\prolem{\ref{matr1}}{ The first statement follows from the definition of out forest. Let us prove the second
statement. Suppose that $(s,i)\in E(\G)$ and $\FF_{k,i}^{j\to s}\ne\emptyset.$ Then $s$ is unreachable from
$i$ in every forest $F_{k}\in \FF_{k,i}^{j\to s}$. After the addition of $(s,i)$ to $F_{k}$, we obtain a
diverging forest with $k+1$ arcs where $i$ is reachable from $j$, i.e., $F_{k}+(s,i)\in\FF_{k+1}^{j\to i}$.

Suppose now that $F_{k+1}\in\FF_{k+1}^{j\to i}$. Let $(p,i)$ be the unique arc directed to $i$ in $F_{k+1}$.
Since $F_{k+1}-(p,i)\in \FF_{k,i}^{j\to p}$, we have $F_{k+1}\in\FF_{k,i}^{j\to p}+(p,i)
\subseteq\cupo_{(s,i)\in E(\G)}(\FF_{k,i}^{j\to s}+(s,i)).$ }

{\bf Corollary from Lemma~\ref{matr1}.} For any digraph and all $i,j=\1n$ and $k=0\cdc n-v,$
$$
q^{k+1}_{ij}=\e(\FF_{k+1}^{j\to i})= \e\Big(\cupo_{(s,i)\in E(\G)} (\FF_{k,i}^{j\to s}+
(s,i))\Big)=\sum_{(s,i)\in E(\G)} \e(\FF_{k,i}^{j\to s}+ (s,i))=\sum_{(s,i)\in E(\G)} \e_{si} q^{k}_{sj,i}.
$$

We now continue proving Proposition~\ref{pro.allk}. Let $LQ_{k}=(a^{k}_{ij}).$ Then for all $i=\1n,$
\begin{eqnarray}
\label{aii.lj} a^{k}_{ii} &\!=&\! \sum_{s=1}^n \l_{is} q^{k}_{si}= \sum_{s \ne i} \l_{is} q^{k}_{si}+\l_{ii}
q^{k}_{ii}= \sum_{s \ne i} \l_{is} q^{k}_{si}-\sum_{s \ne i} \l_{is} q^{k}_{ii}=
\sum_{s=1}^{n} \e_{si} (q^{k}_{ii}-q^{k}_{si})\\
\nonumber &\!=&\! \sum_{s=1}^{n} \e_{si} (q^{k}_{\bar s i}+q^{k}_{si}-q^{k}_{si})= \sum_{s=1}^{n} \e_{si}
q^{k}_{\bar si}=\si\_{k+1}-q^{k+1}_{ii} \ge 0.
\end{eqnarray}

The statement of Proposition~\ref{pro.allk} with respect to the diagonal entries of $L Q_k$ is proved.

For $j\ne i$ we have
\begin{eqnarray}
a^{k}_{ij} &=&\sum_{s=1}^n \l_{is} q^{k}_{sj}= \sum_{s \ne i} \l_{is} q^{k}_{sj}-\sum_{s \ne i} \l_{is}
q^{k}_{ij}= \sum_{s=1}^n \e_{si} q^{k}_{ij}-\sum_{s=1}^n \e_{si}  q^{k}_{sj} \nonumber
\\
&=&\sum_{s=1}^n \e_{si} q^{k}_{ij}- \sum_{s=1}^n \e_{si} (q^{k}_{sj,{\bar i}}+q^{k}_{sj,i})= \sum_{s=1}^n
\e_{si} q^{k}_{ij}- \sum_{s=1}^n \e_{si} q^{k}_{sj,{\bar i}}- \sum_{s=1}^n \e_{si}q^{k}_{sj,i} \nonumber
\\
&=&\e(\GG_1)-\e(\GG_2)-\sum_{s=1}^n \e_{si}q^{k}_{sj,i}, \label{aij1.lj}
\end{eqnarray}
where $\e(\GG_1)=\sum_{s=1}^n\e_{si} q^{k}_{ij},$ $\e(\GG_2)=\sum_{s=1}^n\e_{si} q^{k}_{sj,{\bar i}}$,
$\GG_1$ is the multiset of 2-subgraphs obtained by the addition of all arcs $(s,i)\in E(\G)${\x} to all
possible forests in $\FF_{k}^{j\to i}$, and{\x} $\GG_2$ is the multiset of 2-subgraphs obtained by the
addition of all arcs $(s,i)\in E(\G)${\x} of $\G$ to all possible forests in $\FF_{k,{\bar i}}^{j\to s}$. The
multiset $\GG_1$ consists of the pairs $(H,n_1(H))$ where $n_1(H)\ge1$ is the multiplicity of $H$ in~$\GG_1$;
$\GG_2$ consists of the pairs $(H,n_2(H))$ where $n_2(H)$     is the multiplicity of $H$ in~$\GG_2$. The
notation $H\in\GG_1$ means here $n_1(H)>0$.

Prove that $\GG_1=\GG_2$. In every $H\in\GG_1,$ two arcs are directed to $i$: $(s,i)$ and $(p,i)$ such that
$(p,i)\in E(F_{k}^{j\to i})$, since $j\ne i$. Consider the digraph $H-(p,i)$. It is a forest where $i$ is not
a root and $j$ is the root of the tree that contains~$p$. Then $H-(p,i)\in\FF_{k,{\bar i}}^{j\to p}$ and
$H\in\FF_{k,{\bar i}}^{j\to p}+(p,i)\subseteq\GG_2$. Conversely, suppose that $H=F_{k,{\bar i}}^{j\to
s}+(s,i)\in\GG_2$. Consider $H-(p,i)$ such that $(p,i)\in E(F_{k,{\bar i}}^{j\to s})$. We have
$H-(p,i)\in\FF_{k}^{j\to i}$, hence, $H=F_{k}^{j\to i}+(p,i)\in\GG_1$. For the multisets $\GG_1$ and $\GG_2$,
the multiplicities of all elements, $n_1(H)$ and $n_2(H)$, do not exceeded two. Otherwise, $H-(s,i)$ would
not be a forest. Prove that for any~$H$, $n_1(H)=2$ iff $n_2(H)=2$. Suppose that $H\in\GG_1$ and $n_1(H)=2$.
Then there are vertices $s_1,s_2\in V(\G)$ such that $s_1\ne s_2$, $E(H)$ contains $(s_1,i)$ and $(s_2,i)$,
and $\{H-(s_1,i),\;H-(s_2,i)\}\subseteq \FF_{k}^{j\to i}.$ Since $s_1$ is reachable from $j$ in $H-(s_2,i)$,
$s_1$ is reachable from $j$ in $H$ too. Therefore, $s_1$ is reachable from $j$ in $H-(s_1,i)$ as well.
Similarly, $s_2$ is reachable from $j$ in $H-(s_2,i)$. We obtain $H-(s_1,i)\subseteq\FF_{k,{\bar i}}^{j\to
s_1}$ and $H-(s_2,i)\subseteq\FF_{k,{\bar i}}^{j\to s_2}$, consequently, $n_2(H)=2.$

Suppose that $H\in\GG_2$ and $n_2(H)=2$. Then for some distinct vertices $s_1$ and $s_2,$ there exist
$F_{k,{\bar i}}^{j\to s_1}\in\FF_{k,{\bar i}}^{j\to s_1}$ and $F_{k,{\bar i}}^{j\to s_2}\in\FF_{k,{\bar
i}}^{j\to s_2}$ such that $H=F_{k,{\bar i}}^{j\to s_1}+(s_1,i)=F_{k,{\bar i}}^{j\to s_2}+(s_2,i).$ Note that
(A)~$F_{k,{\bar i}}^{j\to s_1}$ contains the arc $(s_2,i)$, therefore, $s_2$ is not reachable from $i$ in
this forest. Then (B)~$s_2$ is reachable from $j$ in $F_{k,{\bar i}}^{j\to s_1}$, since otherwise $s_2$ would
be reachable from $j$ in $H=F_{k,{\bar i}}^{j\to s_1}+(s_1,i)$, which is false. By (A) and (B), $F_{k,{\bar
i}}^{j\to s_1}\in\FF_{k}^{j \to i}$, and similarly, $F_{k,{\bar i}}^{j\to s_2}\in\FF_{k}^{j \to i}$. Hence,
$n_1(H)=2$.

By (\ref{aij1.lj}) and Corollary from Lemma~\ref{matr1}, at $j\ne i$ we have
\begin{equation}
\label{aij2.lj} a^{k}_{ij}=-\sum_{s\ne i}\e_{si}q^{k}_{sj,i}=-q^{k+1}_{ij}\le 0.
\end{equation}

Proposition~\ref{pro.allk} is proved.
}  %

\prothe{\ref{teo.allk}}{ Let us represent (\ref{system1}) as follows:
\begin{equation}
\label{system2} \cases{ Q_{0}=\si\_0 I=I,                          \cr Q_{1}=\si\_1 I-L Q_{0},
\cr \ldots\ldots\ldots\ldots\ldots\ldots\ldots \cr Q_{k}=\si\_{k}I-LQ_{k-1},                  \cr
\ldots\ldots\ldots\ldots\ldots\ldots\ldots \cr Q_{n-v}=\si\_{n-v}I-LQ_{n-v-1}.            \cr }
\end{equation}

By substituting each equation in the subsequent one, for all $k=2\cdc n-v$ we have
\begin{eqnarray*}
Q_{k} &=& \si\_{k}I
-L\Big(\si\_{k-1}I-\ldots-L(\si\_1 I-L\si\_0 I)\ldots\Big)\\
&=& \si\_{k}I-\si\_{k-1}L+\ldots+{(-1)}^{k}\si\_0 L^k =\sum_{i=0}^{k}\si\_{k-i}{(-L)}^{i}.
\end{eqnarray*}

Theorem~\ref{teo.allk} is proved. }

\prolem{\ref{sumstr}}{ The $i$th row sum of $LQ_k=(a_{ij}^k)$ is
\[
\sum_{j=1}^n a^{k}_{ij}= \sum_{j=1}^n       \sum_{s=1}^n\l_{is}q^{k}_{sj}= \sum_{s=1}^n\l_{is}\sum_{j=1}^n
q^{k}_{sj}= \sum_{s=1}^n\l_{is}\si\_{k}=0.
\]
}

{\bf Proof of Corollary~1 from Theorem~\ref{teo.allk}.} Multiplying both sides of (\ref{psibek}) by any
matrix that commutes with $L$ and using distributivity and associativity of matrix operations,{\x} we get the
required statement.

{\bf Proof of Corollary~2 from Theorem~\ref{teo.allk}.} Consider (\ref{aii.lj}) at $k=n-v$. By virtue of
Proposition~\ref{prop2}, for every $(s,i)\in E(\G)$, $\,\e_{si} q^{n-v}_{\bar si}=0$ holds, i.e.,
$a^{n-v}_{ii}=0$ for all $i\in\{\1n\}$. In this way, Corollary~2 is derived from $a^{n-v}_{ii}=0$,
inequality~(\ref{aij2.lj}), Lemma~\ref{sumstr}, and Corollary~1 from Theorem~\ref{teo.allk}.

{\bf Proof of Corollary~3 from Theorem~\ref{teo.allk}.} Postmultiplying both sides of (\ref{psibek}) by $\q$
and using $L\q=0$ (Corollary~2 from Theorem~\ref{teo.allk}) provides $Q_k\q=\si\_k\q$. By commutativity,
$\q\q_k=\q_k\q=\q$ holds. Using Theorem~4, we also get $Q(\tau)\q=\q Q(\tau)=\q$ for any $\tau>0.$

\prothe{\ref{sumsum7}}{ Substituting~(\ref{psibek}) in~(\ref{iden1}) gives
\begin{eqnarray*}
s\,(I+L)^{-1}&=&\suml_{k=0}^{n-v}\suml_{i=0}^{k}\si\_{k-i}(-L)^i
              = \suml_{i=0}^{n-v}\suml_{k=i}^{n-v}\si\_{k-i}(-L)^i\cr
             &=&\suml_{i=0}^{n-v}\suml_{j=0}^{n-v-i}\si\_{j}(-L)^i
              = \suml_{i=0}^{n-v}s\_{n-v-i}(-L)^i.
\end{eqnarray*}
}

{\bf Proof of Corollary from Theorem~\ref{sumsum7}.} Observe that the digraph resulting from $\G$ by
multiplying the weights of all arcs by $\tau$ has the Kirchhoff matrix $\tau L$, and its total weight of out
forests with $j$ arcs is $\tau^j\si\_j.$ Hence, the required statement follows from Theorem~\ref{sumsum7}.

\prothe{\ref{fulllj}}{ We first prove the following lemma.

\begin{lemma}
\label{lem1} \label{invert.Lj} If a maximum out forest of a digraph $\G$ is a tree$,$ i.e.$,$ the out{\x}
forest dimension of $\G$ is~$1,$ then $L+\vj$ is nonsingular.
\end{lemma}

\prolem{\ref{invert.Lj}}{ Assume, on the contrary, that $\det(L+\vj)=0$. Then there exists a vector ${\bf
b}=(b\_1\cdc b\_n)^{\intercal}\ne{\bf 0}$ such that $(L+\vj)^{\intercal}{\bf b}={\bf 0},$ where ${\bf
0}=(0\cdc 0)^{\intercal}$.

Since, for every digraph of out{\x} forest dimension~1, every column of $\vj$ consists of equal entries
(item~3 of Theorem~\ref{t1.che.ra1}), we obtain
\begin{eqnarray}
\label{lj1} (L+\vj)^{\intercal} {\bf b} &=& \left\|\matrix {
 b_1\l_{11}+\ldots+b_n\l_{n1}+b_1\q_{11}+\ldots+b_n \q_{n1} \cr
 b_1\l_{12}+\ldots+b_n\l_{n2}+b_1\q_{12}+\ldots+b_n \q_{n2} \cr
\ldots\ldots\ldots\ldots\ldots\ldots\ldots\ldots\ldots\ldots\ldots\ldots\ldots\cr
 b_1\l_{1n}+\ldots+b_n\l_{nn}+b_1\q_{1n}+\ldots+b_n \q_{nn} \cr
}\right\|\nonumber\\
&=& \left\|\matrix {b_1\l_{11}+\ldots+b_n\l_{n1}+\q_{11}(b_1+\ldots+b_n)\cr
 b_1\l_{12}+\ldots+b_n\l_{n2}+\q_{12}(b_1+\ldots+b_n)\cr
\ldots\ldots\ldots\ldots\ldots\ldots\ldots\ldots\ldots\ldots\ldots\ldots\cr
 b_1\l_{1n}+\ldots+b_n\l_{nn}+\q_{1n}(b_1+\ldots+b_n)\cr
}\right\|=\left\|\matrix {0\cr 0\cr\vdots\cr 0\cr }\right\|.
\end{eqnarray}

Adding up the components of the vectors that form the last equality and using the identities
$\suml^{n}_{k=1}\!\q_{1k}=1$ (item~1 of Theorem~\ref{t1.che.ra1}) and $\suml^{n}_{k=1}\!\l_{ik}=\!0,$
$i=\1n,$ we deduce $b_1+\ldots+b_n=0$. Replacing the last equation of~(\ref{lj1}) with this equality, we
obtain the following system of equations{\x} in $b_1\cdc b_n$:
\begin{equation}
\label{b1bn} \cases{
 \: b_1\l_{11}+\ldots+b_n \l_{n1}    &$\!\!\!\!=0$ \cr
 \: b_1\l_{12}+\ldots+b_n \l_{n2}    &$\!\!\!\!=0$ \cr
 \: \ldots\ldots\ldots\ldots\ldots\ldots\ldots     \cr
    b_1\l_{1n-1}+\ldots+b_n \l_{nn-1}&$\!\!\!\!=0$ \cr
    b_1+\ldots+b_n                   &$\!\!\!\!=0$.\cr
}
\end{equation}

Let $M$ be the transposed matrix of coefficients of the system~(\ref{b1bn}). Expand the determinant of $M$ in
the last column (which consists of ones). By the matrix-tree theorem for digraphs (see,
e.g.,~\cite{HarNoCa,Harary}), if there exists a spanning tree diverging from $i$, then the cofactor of
$\l_{in}$ in $L$---which equals the cofactor of the $i$th entry of the last column of $M$---is positive,
whereas, in the opposite case, it is zero. By the hypothesis of this lemma, the maximum out forests of $\G$
are diverging trees. Therefore, the cofactor of at least one entry of the last column of $M$ is positive.
Hence, $\det M$ expanded as above is also positive. Thereby, ${\rm rank}\,M=n$, whence the unique solution of
the system~(\ref{b1bn}) is $b_1=\cdots=b_n=0$. This contradiction proves Lemma~\ref{invert.Lj}. }

Suppose that $V_1\cdc V_p$ are the vertex sets of the weak components of~$\G$. Without loss of generality, we
assume that the vertices of $V_1$ are numbered first, the vertices of $V_2$ are numbered next, etc. (at any
other numeration, the corresponding permutation of the rows and columns preserves the rank of $L+\vj$). By
virtue of Theorem~$2'$ in \cite{Che.ra1}, $L+\q$ is a block diagonal matrix with $p$ blocks, which is
expressed in the following manner~\cite{HoJo}:
\begin{equation}
\label{ass1} L+\q=\oplus\suml^{p}_{s=1}A_{ss},
\end{equation}
where submatrix $A_{ss}$ corresponds to the $s$th weak component. A block diagonal matrix $L+\q$ is
nonsingular if and only if every $A_{ss}$ is nonsingular, moreover,
\begin{equation}
\label{ass2} {\rm det}(L+\q)=\prod^{p}_{s=1}{\rm det}A_{ss}.
\end{equation}

Let $\G_s$ be the restriction of $\G$ to $V_s,$ $s\in \{1\cdc p\}$; $L^{(s)}, Q^{(s)}_{n-v}$, and $\vj^{(s)}$
will be the matrices constructed for~$\G_s$. Then
\begin{equation}
\label{ass3} A_{ss}=L^{(s)}+\vj^{(s)}.
\end{equation}
Indeed, $L^{(s)}$ coincides with the $s$th block of~$L$, whereas, by item~5 of Theorem~\ref{t1.che.ra1}, the
$s$th block of $Q_{n-v}$ is proportional to $Q^{(s)}_{n-v}$ (the proportionality factor being the total
weight of all forests in $\G_{-V_s}$). Consequently, the $s$th block of $\vj$ coincides with~$\vj^{(s)}$.

By virtue of (\ref{ass2}) and (\ref{ass3}), it is sufficient to prove the statement of the theorem for the
case of $p=1$. Let $\G$ consist of a single weak component. Suppose, without loss of generality, that $V(\G)$
is indexed in such a way that the first numbers are attached to the vertices in $K_1$, the subsequent numbers
to the vertices in $K_2$, etc. The last numbers are given to the vertices in $V(\G)\backslash\ktil$. By
item~2 of Theorem~\ref{t1.che.ra1}, $L+\vj$ is a block lower triangular matrix with $v+1$ blocks: $K_1$
corresponds to the first block, $K_v$ corresponds to the $v$th block; $V(\G)\backslash\ktil$ corresponds to
the last block.

The determinant of the block triangular matrix $L+\vj$ is the product of the determinants of its diagonal
blocks, whereas its rank is no less than the sum of ranks of the diagonal blocks (see, e.g.,~\cite{HoJo}).
Note that, as well as in the case of weak components of $\G$ (see above), the block of $L+\vj$ corresponding
to an undominated knot $K_i$, coincides with the matrix $L(\G_{K_i})+\vj(\G_{K_i})$ constructed for
$\G_{K_i}$. To demonstrate this, it suffices to use item~3 of Theorem~\ref{t1.che.ra1}. Every maximum out
forest of an undominated knot is a diverging tree (out arborescence). Therefore, the nonsingularity of the
diagonal blocks of $L+\vj$ that correspond to the undominated knots follows from Lemma~\ref{invert.Lj}.

The $(v+1)$st{\x} block of $L+\vj$ coincides with the corresponding block of $L$, since the last block of
$\vj$ is zero (by item~2 of Theorem~\ref{t1.che.ra1}). The last block of $L$ is nonsingular. Indeed, by
Theorem~6 in~\cite{Fiedler}, its determinant is the weight of the set of out forests in $\G$ where $\ktil$ is
the set of roots. This weight is strictly positive, since the indicated set of forests is nonempty: such
forests can be obtained from the maximum out forests of $\G$ by the removal of all arcs between the vertices
within~$\ktil$.

Thus, the diagonal blocks of $L+\vj$ are nonsingular, hence, $L+\vj$ is nonsingular. The theorem is proved. }

{\bf Proof of Corollary from Theorem~\ref{fulllj}.} First, we prove the following lemma.

\begin{lemma}
\label{invert.Laj} For every weighted digraph $\G$ of out{\x} forest dimension $1$ and any $\aa\neq 0,$ the
matrix $L+\aa\vj$ is nonsingular.
\end{lemma}

\prolem{\ref{invert.Laj}}{ This lemma is proved by the same argument as Lemma~\ref{invert.Lj} with the only
difference that the analogue of (\ref{lj1}) takes here a more general form:
\[
(L+\aa\vj)^{\intercal}{\bf b} =\left\|\matrix{ b_1\l_{11}+\ldots+b_n\l_{n1}+\aa\q_{11}(b_1+\ldots+b_n)\cr
b_1\l_{12}+\ldots+b_n\l_{n2}+\aa\q_{12}(b_1+\ldots+b_n)\cr
\ldots\ldots\ldots\ldots\ldots\ldots\ldots\ldots\ldots\ldots\ldots\ldots\cr
b_1\l_{1n}+\ldots+b_n\l_{nn}+\aa\q_{1n}(b_1+\ldots+b_n)\cr }\right\|=\left\|\matrix {0\cr 0\cr\vdots\cr 0\cr
}\right\|.
\]
}

To complete the proof of Corollary from Theorem~\ref{fulllj}, note that for any $\aa\ne0$, the matrix
$L+\aa\vj$, as well as $L+\vj$, is a block lower triangular matrix with $v+1$ blocks. By item~2 of
Theorem~\ref{t1.che.ra1}, its $(v+1)$st{\x} diagonal block coincides with the corresponding block of $L+\vj$.
Using Lemma~\ref{invert.Laj}, we conclude that the other diagonal blocks are also nonsingular.

\prothe{\ref{groupinv}}{ By Theorems~\ref{t2.che.ra1} and~\ref{t3.che.ra1}, $(L+\oj)\oj=\oj$. Premultiplying
both sides of this identity by {\x}${(L+\oj)}^{-1}$ (which exists by Theorem~\ref{fulllj}), we have
\begin{equation}
\label{vsp1} \oj={(L+\oj)}^{-1}\oj.
\end{equation}

Similarly,
\begin{equation}
\label{vsp1a} \q(L+\oj)^{-1}=\oj.
\end{equation}

Denote{\x} $\qtil=(L+\vj)^{-1}-\vj$. Using (\ref{vsp1}) and Theorems~\ref{t2.che.ra1} and~\ref{t3.che.ra1},
we obtain
\begin{equation}
\label{vsp2} \qtil L={(L+\oj)}^{-1}L-\oj L
       ={(L+\oj)}^{-1}(L+\oj-\oj)
       =I-{(L+\oj)}^{-1}\oj
       =I-\oj.
\end{equation}

Similarly,
\begin{equation}
\label{vsp3} L\qtil=I-\oj.
\end{equation}
By (\ref{vsp1}) and Theorem~\ref{t2.che.ra1},
\[
\qtil\oj={(L+\oj)}^{-1}\oj-\oj^2=0.
\]
Consequently, for any $\alpha\ne0,$
\[
(\qtil+\alpha^{-1}\oj)(L+\alpha\oj)=I-\oj+\oj=I.
\]
Hence, $\qtil+\alpha^{-1}\oj=(L+\alpha\oj)^{-1}$ and
\begin{equation}
\label{vsp4} \qtil=(L+\aa\vj)^{-1}-\aa^{-1}\vj
\end{equation}
for any $\alpha\ne0$.

By (\ref{vsp2}) and (\ref{vsp3}), $L\qtil=\qtil L=I-\oj$, thus $\qtil=(L+\aa\vj)^{-1}-\aa^{-1}\vj$ satisfies
condition~(5) in the definition of group inverse. Let us prove that conditions~(1) and~(2) (common with the
definition of Moore-Penrose inverse) are also fulfilled. Making use of Theorems~\ref{t2.che.ra1}
and~\ref{t3.che.ra1} and identities
(\ref{vsp1a})--(\ref{vsp4}), we obtain
\begin{eqnarray*}
L\qtil L     &=& L(I-\oj)=L,\\
\qtil L\qtil &=&(I-\oj)\qtil
              =\qtil-\oj \qtil
              =\qtil-\oj{(L+\oj)}^{-1}+\oj^2
              =\qtil-\oj+\oj=\qtil,
\end{eqnarray*}
which completes the proof. }

\propro{\ref{qtiLim}}{ Using Theorems~\ref{t2.che.ra1},~\ref{t3.che.ra1}, and~\ref{t5.che.le2} and identities
$(I+\tau L)^{-1}\q=\q$ (Corollary~3 from Theorem~\ref{teo.allk}) and~(\ref{vsp1}), we obtain
\begin{eqnarray*}
&&\left(\liml_{\tau\to\infty}\tau\left({(I+\tau L)}^{-1}-
   \oj\right)+\oj\right)(L+\oj)\\
&&=\liml_{\tau\to\infty}
   \tau\left({(I+\tau L)}^{-1}L+{(I+\tau L)}^{-1}\oj-\oj L-
   {\oj}^2\right)+\oj L+\oj^2\\
&&=\liml_{\tau\to\infty}\tau{(I+\tau L)}^{-1}L+\oj
  =\liml_{\tau\to\infty}{(I+\tau L)}^{-1}(I+\tau L-I)+\oj\\
&&=I-\liml_{\tau\to\infty}{(I+\tau L)}^{-1}+\oj=I.
\end{eqnarray*}

Postmultiplying the first and last expressions by $(L+\q)^{-1}$ and using Theorem~\ref{groupinv} we obtain
the required equation. }

\propro{\ref{topol}}{ Using Proposition~\ref{qtiLim} and Theorem~4, for any $i,j\in V(\G),$ find the limit
\begin{eqnarray*}
\lpr_{ij} &=&\liml_{\tau\to\infty}\tau\Big(q\_{ij}(\tau)-\vj_{ij}\Big)=
   \liml_{\tau\to\infty}\tau\left(\frac{\suml^{n-v}_{k=0}\tau^k{\e}
(\FF^{j\to i}_k)} {\suml^{n-v}_{k=0}\tau^k\si\_k}-\vj_{ij}\right)
\\
&=&\liml_{\tau\to\infty}\frac{\suml^{n-v}_{k=0}\tau^{k+1}{\e} (\FF^{j\to
i}_k)-\suml^{n-v}_{k=0}\tau^{k+1}\si\_k\vj_{ij}} {\suml^{n-v}_{k=0}\tau^k\si\_k}
\\
&=&\liml_{\tau\to\infty}\frac{\suml^{n-v-1}_{k=0}\tau^{k+1}{\e} (\FF^{j\to i}_k) +\tau^{n-v+1}{\e}(\FF^{j\to
i}_{n-v}) -\suml^{n-v-1}_{k=0}\tau^{k+1}\si\_k\vj_{ij} -\tau^{n-v+1}\si\_{n-v}\vj_{ij}}
 {\suml^{n-v}_{k=0}\tau^k\si\_k}.
\end{eqnarray*}

By the definition of $\vj$, we have $\tau^{n-v+1}{\e}(\FF^{j \to i}_{n-v})-
                                    \tau^{n-v+1}\si\_{n-v}\vj_{ij}=0$.
Therefore, in view of $\si\_{n-v}\ne0$, this yields{\x}
\[
\lpr_{ij}=\frac{\e(\FF^{j\to i}_{n-v-1})-\si\_{n-v-1}\vj_{ij}}{\si\_{n-v}},
\]
which completes the proof of Proposition~\ref{topol}.}

\prolem{\ref{commut}} {By virtue of the identity $\vj L=0$ (Theorem~\ref{t3.che.ra1}), matrices
$LL^{\intercal}$ and $ZZ^{\intercal}=\vj^{\intercal}\!\vj+L L^{\intercal}$ commute, i.e., $L L^{\intercal}
(\vj^{\intercal}\!\vj+L L^{\intercal})=(\vj^{\intercal}\!\vj+L L^{\intercal})L
L^{\intercal}={(LL^{\intercal})}^2$. Premultiplying and postmultiplying both sides of the first equality by
${(ZZ^{\intercal})}^{-1}=(\vj^{\intercal}\!\vj+LL^{\intercal})^{-1}$, we obtain the desired
$(\vj^{\intercal}\!\vj +LL^{\intercal})^{-1}LL^{\intercal}=LL^{\intercal}(\vj^{\intercal}\!\vj
+LL^{\intercal})^{-1}$. The second statement is proved similarly. }

\prothe{\ref{pseudoor}}{ Let $\ltil=L^{\intercal}{(ZZ^{\intercal})}^{-1}$. Prove that the following four
conditions are satisfied:

(1) $L\ltil L=L$;

(2) $\ltil L \ltil=\ltil$;

(3) $L\ltil$ is symmetric;

(4) $\ltil L$ is symmetric.

Condition 1. Using Lemma~\ref{commut} and the identity $\vj L=0$ (Theorem~\ref{t3.che.ra1}), we obtain
\begin{eqnarray*}
L \ltil L &=&LL^{\intercal}(\vj^{\intercal}\!\vj+L L^{\intercal})^{-1}L
=((LL^{\intercal}+\vj^{\intercal}\!\vj) -\vj^{\intercal}\!\vj)(\vj^{\intercal}\!\vj +LL^{\intercal})^{-1}L
\\
&=&(I-\vj^{\intercal}\!\vj (\vj^{\intercal}\!\vj+LL^{\intercal})^{-1})L
=L-(\vj^{\intercal}\!\vj+LL^{\intercal})^{-1}\vj^{\intercal}\!\vj L=L.
\end{eqnarray*}

Condition 2. Using Lemma~\ref{commut} and the identity $L^{\intercal}\!\vj^{\intercal}=0$, we have{\x}
\begin{eqnarray*}
\ltil L\ltil &=&\ltil L L^{\intercal}(\vj^{\intercal}\!\vj+L L^{\intercal})^{-1} =\ltil((L
L^{\intercal}+\vj^{\intercal}\!\vj) -\vj^{\intercal}\!\vj)(\vj^{\intercal}\!\vj+L L^{\intercal})^{-1}
\\
&=&\ltil(I-\vj^{\intercal}\!\vj(\vj^{\intercal}\!\vj+L L^{\intercal})^{-1})
=\ltil-L^{\intercal}(\vj^{\intercal}\!\vj +L L^{\intercal})^{-1}\vj^{\intercal}\!\vj(\vj^{\intercal}\!\vj +L
L^{\intercal})^{-1}
\\
&=&\ltil-L^{\intercal}\!\vj^{\intercal}\!\vj(\vj^{\intercal}\!\vj +L L^{\intercal})^{-2} =\ltil.
\end{eqnarray*}

Condition 3. By Corollary from Lemma~\ref{commut}, the matrix
$L\ltil=LL^{\intercal}(\vj^{\intercal}\!\vj+LL^{\intercal})^{-1}$ is symmetric.

Condition 4. Since $(\vj^{\intercal}\!\vj+LL^{\intercal})^{-1}$ is symmetric, $\ltil
L=L^{\intercal}(\vj^{\intercal}\!\vj+L L^{\intercal})^{-1}L$ is symmetric as well (see, e.g., \cite{HoJo},
Theorem~4.1.3).

Theorem~\ref{pseudoor} is proved. }

\propro{\ref{annul}}{ By Theorem~\ref{teo.allk}, $Q_{n-v}=\suml_{i=0}^{n-v} \si\_{n-v-i} {(-L)}^{i}$. Using
Corollary~2 from Theorem~\ref{teo.allk}, we obtain
\[
p'_L(L)=L\suml_{i=0}^{n-v}\si\_{n-v-i}{(-L)}^{i}=L Q_{n-v}=0.
\]
}

\propro{\ref{otledostup}}{ Under the reversal of all arcs in a digraph, diverging and converging forests
change places. Therefore, $P\_1(\tau)$ and $P\_2(\tau)$ are dual to each other.

{\it Nonnegativity} and {\it disconnection condition\/} for $P\_1(\tau)$ and $P\_2(\tau)$ and also item~{\rm
1} of {\it diagonal maximality\/} for $P\_1(\tau)$ and item~{\rm 2} of {\it diagonal maximality\/} for
$P\_2(\tau)$ are proved with the help of Theorem~4 by the same argument as the corresponding conditions
in~\cite{CheSha97}.

For any $i,k,t\in V(G),$ if $\G$ contains a path from $i$ to $k,$ $i\ne k\ne t,$ and every path from $i$ to
$t$ includes $k,$ then $\FF^{i\to t}\subset\FF^{i\to k}$. This inclusion and Theorem~4 imply the fulfillment
of item~1 of {\it transit property\/} for $P\_1(\tau)$. Item~2 of {\it transit property\/} for $P\_2(\tau)$
is proved similarly.

It was established in \cite{Che.ra1} (in the proof of Proposition~17) that if the weight of some arc $(k,t)$
is increased by $\D\e_{kt}$ and the weights of all other arcs are preserved, then the increments of the
entries of the matrix $P\_1(\tau)=\left(p^{(1)}_{ij}(\tau)\right)$ are expressed as follows:
\begin{equation}
\label{deltaptau2} \D\,p^{(1)}_{ij}(\tau)
=\frac{p^{(1)}_{tj}(\tau)\left(p^{(1)}_{ik}(\tau)-p^{(1)}_{it}(\tau)\right)}
      {(\D\e\_{kt}\tau)^{-1}+p^{(1)}_{tt}(\tau)-p^{(1)}_{tk}(\tau)},
\quad i,j=\1n.
\end{equation}

By {\it reversal property}, the corresponding formula for $P\_2(\tau)=\left(p^{(2)}_{ij}(\tau)\right)$ is
\begin{equation}
\label{delta1210} \D\,p^{(2)}_{ij}(\tau)
=\frac{p^{(2)}_{ik}(\tau)\left(p^{(2)}_{tj}(\tau)-p^{(2)}_{kj}(\tau)\right)}
      {(\D\e\_{kt}\tau)^{-1}+p^{(2)}_{kk}(\tau)-p^{(2)}_{tk}(\tau)},
\quad i,j=\1n,
\end{equation}
where $p^{(2)}_{ij}(\tau)$ is the total weight of in forests wherein $j$ belongs to a tree converging to~$i$.

Transcribe (\ref{deltaptau2}) for $\D\,p^{(1)}_{kt}(\tau)$:
\[
\D\,p^{(1)}_{kt}(\tau) =\frac{p^{(1)}_{tt}(\tau)\left(p^{(1)}_{kk}(\tau)-p^{(1)}_{kt}(\tau)\right)}
      {(\D\e\_{kt}\tau)^{-1}+p^{(1)}_{tt}(\tau)-p^{(1)}_{tk}(\tau)}.
\]
Since $p^{(1)}_{ii}(\tau)>p^{(1)}_{ij}(\tau)$ for all $i,j=\1n$, we deduce $\D\,p^{(1)}_{kt}(\tau)>0$, i.e.,
$P\_1(\tau)$ satisfies the first part of item~1 of {\it monotonicity}.

Now represent $\D\,p^{(2)}_{kt}(\tau)$ using (\ref{delta1210}):
\[
\D\,p^{(2)}_{kt}(\tau) =\frac{p^{(2)}_{kk}(\tau)\left(p^{(2)}_{tt}(\tau)-p^{(2)}_{kt}(\tau)\right)}
      {(\D\e\_{kt}\tau)^{-1}+p^{(2)}_{kk} (\tau)-p^{(2)}_{tk}(\tau)}.
\]
Since $p^{(2)}_{ii}(\tau)>p^{(2)}_{ji}(\tau)$ for all $i,j=\1n$, we have $\D\,p^{(2)}_{kt}(\tau)>0$, i.e.,
$P\_2(\tau)$ also satisfies the first part of item~1 of {\it monotonicity}.

Transcribing (\ref{deltaptau2}) for $\D\,p^{(1)}_{kt}(\tau)$, $\D\,p^{(1)}_{ki}(\tau)$, and
$\D\,p^{(1)}_{ti}(\tau)$ and using item~1 of diagonal maximality, we conclude that $P\_1(\tau)$ satisfies
item~2 of monotonicity. Comparing the expressions for $\D\,p^{(1)}_{it}(\tau)$ and $\D\,p^{(1)}_{ik}(\tau)$
and using item~1 of transit property, we obtain item~$3b$ of monotonicity. Along the same lines, the
statement of item~$3a$ holds if and only if $\Big(p^{(1)}_{kk}(\tau)-p^{(1)}_{kt}(\tau)\Big)-
 \Big(p^{(1)}_{ik}(\tau)-p^{(1)}_{it}(\tau)\Big)>0$.
Since triangle inequality for accessibility measures is violated by $P\_1(\tau)$ (as well as by all
asymmetric measures), item~$3a$ of monotonicity is not true. Similarly, $P\_2(\tau)$ satisfies items~3
and~$2b$, but violates item~$2a$ of monotonicity.

Consider the digraph with vertex set $\{i,j,k,t\}$, arc set $\{(i,k),(k,t),(t,j)\}$, and arc weights
$\e(i,k)=4$, $\e(k,t)=1$, and $\e(t,j)=4$. Then $P\_1(1)$ is given by
\begin{equation}
P\_1(1)=\begin{array}{l@{\:}llll@{\:}l@{\;}l}
\smallskip
    & i & j    & k   & t&\\
\moo& 1 & 0.32 & 0.8 & 0.4&\mo&i\\
\moo& 0 & 0.2  & 0   & 0 &\mo&j\\
\moo& 0 & 0.08 & 0.2 & 0.1&\mo&k\\
\moo& 0 & 0.4  & 0   & 0.5&\mo&t\\
{}\\
\end{array}.
\end{equation}

Since $p^{(1)}_{it}(1)>p^{(1)}_{kt}(1),$ item~2 of transit property is not satisfied.

Compare the increments $\D\,p^{(1)}_{kt}(1)$ and $\D\,p^{(1)}_{ij}(1)$ for an arbitrary $\D\e_{kt}>0${\rm:}
\[
\D\,p^{(1)}_{kt} =\frac{p^{(1)}_{tt}(1)\left(p^{(1)}_{kk}(1)-p^{(1)}_{kt}(1)\right)}
{(\D\e\_{kt})^{-1}+p^{(1)}_{tt}(1)-p^{(1)}_{tk}(1)}, \quad\quad \D\,p^{(1)}_{ij}
=\frac{p^{(1)}_{tj}(1)\left(p^{(1)}_{ik}(1)-p^{(1)}_{it}(1)\right)}
{(\D\e\_{kt})^{-1}+p^{(1)}_{tt}(1)-p^{(1)}_{tk}(1)}.
\]

Since $p^{(1)}_{tt}(1)\left(p^{(1)}_{kk}(1)-p^{(1)}_{kt}(1)\right)=0.5(0.2-0.1)=0.05$,
$p^{(1)}_{tj}(1)\left(p^{(1)}_{ik}(1)-p^{(1)}_{it}(1)\right)=0.4(0.8-0.4)=0.16$ and the common denominator is
positive, $\D\,p^{(1)}_{kt}(1)<\D\,p^{(1)}_{ij}(1)$ follows,{\x} i.e., the second part of item~1 of
monotonicity is not satisfied.

It is easy to verify that item~1 of transit property and the second part of item~1 of monotonicity are
violated by{\x} $P\_2(\tau)$ in the same example.

Let us show that $P\_3(\tau)$ does not satisfy {\it triangle inequality for accessibility measures}. Consider
the digraph with vertex set $\{i,j,k,t\}$, arc set $\{(i,j),(j,k)$, $(k,t),(t,i)\}$, and arc weights
$\e(i,j)=1,\;\e(j,k)=10,\;\e(k,t)=10,$ and $\e(t,i)=1$. Here, $P\_3(1)$ is as follows:
\begin{equation}
P\_3(1)=\begin{array}{l@{\:}llll@{\:}l@{\;}l}
\smallskip
    & i      & j      & k      & t     &\\
\moo& 0.6302 & 0.2233 & 0.1693 & 0.2233&\mo&i\\
\moo& 0.2233 & 0.3724 & 0.1823 & 0.2747&\mo&j\\
\moo& 0.1693 & 0.1823 & 0.1146 & 0.1823&\mo&k\\
\moo& 0.2233 & 0.2747 & 0.1823 & 0.3724&\mo&t\\
{}\\
\end{array}.
\end{equation}

Triangle inequality for accessibility measures is violated, because
$p_{ki}^{(3)}(1)+p_{kj}^{(3)}(1)-p_{ij}^{(3)}(1)>p_{kk}^{(3)}(1)$. Symmetry is a necessary condition of
triangle inequality for accessibility measures. Therefore, $P\_1(\tau)$ and $P\_2(\tau)$ do not satisfy{\x}
this inequality either. }

\propro{\ref{p1.2}}{ 1.~To construct the required forest, it suffices to take a maximum out forest and to
delete any arc directed to $i$ (if such an arc exists) or {\em any\/} arc in the opposite case.

2.~For the given path (chain subgraph) and any maximum out forest in $\G$, consider their join and remove all
arcs of the out forest that are directed to the vertices of the path but do not belong to the path. The
resulting subgraph contains neither circuits nor vertices with indegree greater than one, i.e., it is an out
forest. It contains at least $n-v-1$ arcs, consequently it belongs to $\FF_{n-v-1}\cup\FF_{n-v}.$}

\propro{\ref{pogushe}}{ Under the reversal of all arcs in $\G$, diverging and converging forests change
places. Hence, $P\_1(\aa)$ and $P\_2(\aa)$ are dual to each other. {\it Disconnection condition\/} follows
from item~2 of Proposition~\ref{p1.2}.

{\it Nonnegativity\/} follows from the fact that the entries of $\q_{n-v}$ and $\q_{n-v-1}$ are proportional
to the weights of some sets and from{\x} the nonnegativity of the coefficients in~(\ref{f21}).

Item~1 of {\it diagonal maximality\/} and item~1 of {\it transit property\/} for $P\_1(\aa)$ in the nonstrict
form (as well as items~2 of these conditions for $P\_2(\aa)$ in the nonstrict form) follow from the nonstrict
inclusion of the sets of forests that determine the entries of matrices $Q_{n-v}$ and $Q_{n-v-1}$ under
comparison.

{\it Diagonal maximality} (in the strict form) is violated, for example, for digraphs with vertex set
$\{j,i,k,t\}$, arc set $\{(j,i),(i,k)(k,t)\}$, and arc weights $\e(j,i)=4,\;\e(i,k)=1,$ and $\e(k,t)=1$. The
matrices of out forests $Q_{n-v}$   and $Q_{n-v-1}$, the matrices of in forests $S_{n-v}$   and $S_{n-v-1}$,
and the matrices $P\_1(\aa)$ and $P\_2(\aa)$ at $\aa=4/13$ are as follows:
\[\hspace{-4.5em}
Q_{n-v}=\begin{array}{l@{\:}llll@{\:}l@{\;}l}
\smallskip
    & j & i & k & t&\\
\moo& 4 & 0 & 0 & 0&\mo&j\\
\moo& 4 & 0 & 0 & 0&\mo&i\\
\moo& 4 & 0 & 0 & 0&\mo&k\\
\moo& 4 & 0 & 0 & 0&\mo&t\\
{}\\
\end{array},
\;\;Q_{n-v-1}=\begin{array}{l@{\:}llll@{\:}l@{\;}l}
\smallskip
    & j & i & k & t&\\
\moo& 9 & 0 & 0 & 0&\mo&j\\
\moo& 8 & 1 & 0 & 0&\mo&i\\
\moo& 4 & 1 & 4 & 0&\mo&k\\
\moo& 0 & 1 & 4 & 4&\mo&t\\
{}\\
\end{array},
\;\;P\_1(\aa)=\begin{array}{l@{\:}l@{\:}l@{\:}l@{\:}l@{\:}l@{\;}l}
\smallskip
    & j    & i    & k    & t   &\\
\moo & 3.25 & 3    & 2    & 1   &\mo & j\\
\moo & 0    & 0.25 & 0.25 & 0.25&\mo & i\\
\moo & 0    & 0    & 1    & 1   &\mo & k\\
\moo & 0    & 0    & 0    & 1   &\mo & t\\
{}\\
\end{array},
\]
\[
\hspace{-4em} S_{n-v}=\begin{array}{l@{\:}llll@{\:}l@{\;}l}
\smallskip
    & j & i & k & t&\\
\moo & 0 & 0 & 0 & 4&\mo &j\\
\moo & 0 & 0 & 0 & 4&\mo &i\\
\moo & 0 & 0 & 0 & 4&\mo &k\\
\moo & 0 & 0 & 0 & 4&\mo &t\\
{}\\
\end{array},
\;\;S_{n-v-1}=\begin{array}{l@{\:}llll@{\:}l@{\;}l}
\smallskip
     & j & i & k & t&\\
\moo & 1 & 4 & 4 & 0&\mo & j\\
\moo & 0 & 4 & 4 & 1&\mo & i\\
\moo & 0 & 0 & 4 & 5&\mo & k\\
\moo & 0 & 0 & 0 & 9&\mo & t\\
{}\\
\end{array},
\;\;P\_2(\aa)=\begin{array}{l@{\:}llll@{\:}l@{\;}l}
\smallskip
    & j    & i & k & t   &\\
\moo & 0.25 & 1 & 1 & 1   &\mo & j\\
\moo & 0    & 1 & 1 & 1.25&\mo & i\\
\moo & 0    & 0 & 1 & 2.25&\mo & k\\
\moo & 0    & 0 & 0 & 3.25&\mo & t\\
{}\\
\end{array}.
\]

In this example, item~1 of transit property for $P\_1(\aa)$ is violated (since
$p^{(1)}_{ik}(\aa)=p^{(1)}_{it}(\aa)$), and  so is item~2 of transit property for $P\_2(\aa)$ (since
$p^{(2)}_{ik}(\aa)=p^{(2)}_{jk}(\aa)$).

Let us demonstrate that {\it triangle inequality for accessibility measures\/} is not satisfied by
$P\_3(\aa)=(P\_1(\aa)+P^{\intercal}_1(\aa)+P\_2(\aa)+P^{\intercal}_2(\aa))/4$ in{\x} this example. Indeed,
\[
P=P\_3(\aa)=\begin{array}{l@{\:}llll@{\:}l@{\;}l}
\smallskip
     & j    & i      & k      & t     &\\
\moo & 1.75 & 1      & 0.75   & 0.5   &\mo & j\\
\moo & 1    & 0.625  & 0.3125 & 0.375 &\mo & i\\
\moo & 0.75 & 0.3125 & 1      & 0.8125&\mo & k\\
\moo & 0.5  & 0.375  & 0.8125 & 2.125 &\mo & t\\
{}\\
\end{array},
\]
and since $p_{ij}^{(3)}(\aa)+p_{it}^{(3)}(\aa)-p_{jt}^{(3)}(\aa) =0.875>p_{ii}^{(3)}(\aa)=0.625$, this
condition is violated.

{\it Monotonicity\/} is not satisfied on undirected graphs (see Proposition~10 in~\cite{CheSha981}),
hence,{\x} it is violated for digraphs also. }

\revred{V.A. Lototskii}

\end{document}